\documentclass{article}
\usepackage{amsfonts, amsmath, graphicx, float}

\newtheorem{thm}{Theorem}
\newtheorem{dfn}{Definition}[subsection]
\newtheorem{prp}[dfn]{Proposition}
\newtheorem{lmm}[dfn]{Lemma}
\newtheorem{cor}[dfn]{Corollary}

\makeatletter

\@addtoreset{equation}{section}
\makeatother

\def \N {{\mathbb N}}

\def \R {{\mathbb R}}

\newcommand{\one}{{\tt 1 \hspace{-0.8ex} \tt l}}
\newcommand{\qed}{\hspace*{2ex} \hfill $\Box$}

\title{
  
  Sampling the Fermi statistics\\
  and other conditional product measures\footnotemark[1]
}
\author{
  
  A. Gaudilli\`ere\footnotemark[2]
  \and
  
  J. Reygner\footnotemark[3]
}
\date{}

\footnotetext[1]{
  The whole research for this paper
  was done at the Dipartimento di matematica
  dell'universit\`a Roma Tre,
  during A.G.'s post-doc
  and J.R.'s research internship. 
  This work was supported 
  by the European Research Council
  through the ``Advanced Grant''  
  PTRELSS 228032.
} 
\footnotetext[2]{
  gaudilli@cmi.univ-mrs.fr
  -- LATP,
  Universit\'e de Provence,
  CNRS,
  39 rue F. Joliot Curie,
  13013 Marseille, France.
}
\footnotetext[3]{
  julien.reygner@polytechnique.edu
  -- CMAP, \'Ecole Polytechnique,
  Route de Saclay,
  91120 Pa\-laiseau, France.
}

\begin{document}

\maketitle

\begin{abstract}
  Through a Metropolis-like algorithm
  with single step computational cost
  of order one, we build a Markov
  chain that relaxes to the canonical
  Fermi statistics for $k$ non-interacting
  particles among $m$ energy levels.
  Uniformly over the temperature
  as well as the energy values and degeneracies
  of the energy levels
  we give an explicit upper bound
  with leading term $km\ln k$ 
  for the mixing time of the dynamics.
  We obtain such construction
  and upper bound as a special
  case of a general result
  on (non-homogeneous) products of
  ultra log-concave measures
  (like binomial or Poisson laws)
  with a global constraint.
  As a consequence
  of this general result
  we also obtain  
  a disorder-independent
  upper bound on the mixing time
  of a simple exclusion process
  on the complete graph with site disorder.
  This general result
  is based on an elementary
  coupling argument,
  illustrated in a simulation appendix
  and extended
  to (non-homogeneous) products
  of log-concave measures.


\end{abstract}

\newpage

\section{From the Fermi statistics
  to general conditional products of log-concave measures
}
\subsection{Sampling the Fermi statistics}

Given two positive integers $k$ and $m$,
given a non-negative real number $\beta$,
given $m$ real numbers $v_1$,~\dots, $v_m$
and given $m$ integers $n_1$,~\dots, $n_m$ such that
\begin{equation}
  n:= \sum_{j = 1}^m n_j \geq k,
\end{equation}
the canonical Fermi statistics 
at {\em inverse temperature} $\beta$
for $k$ non-interacting particles
among the $m$ {\em energy levels} $1$,~\dots, $m$,
with {\em energy values} $v_1$,~\dots, $v_m$
and {\em degeneracies} $n_1$,~\dots, $n_m$ is
the conditional probability measure on
\begin{equation}
  {\cal X}_{k,m} := \left\{
    (k_1,\dots,k_m)\in \N^m :\:
    k_1 + \cdots + k_m = k
  \right\}
\end{equation}
given by 
\begin{equation}
  \nu := \mu(\cdot|{\cal X}_{k,m})
\end{equation}
with $\mu$ the product measure
on $\N^m$ such that
\begin{equation}
  \mu(k_1,\dots, k_m) := 
  \frac{1}{Z}
  \prod_{j=1}^m \binom{n_j}{k_j}\exp\{-\beta k_jv_j\},
  \label{defmu}
\end{equation}
\begin{equation}
  Z := 
  \sum_{k_1,\dots,k_m}
  \prod_{j=1}^m \binom{n_j}{k_j}\exp\{-\beta k_jv_j\},
\end{equation}
with $\binom{n_j}{k_j} = 0$ whenever $k_j > n_j$. 
In other words, $\nu$ is a (non-homogeneous) product
of binomial laws in $k_1$,~\dots, $k_m$
with the global constraint
\begin{equation}
  k_1 + \dots + k_m = k
  \label{canonical}
\end{equation}
and we can write
\begin{equation}
  \nu(k_1,\dots,k_m) = \frac{1}{Q}\prod_{j=1}^m e^{\phi_j(k_j)},
  \quad
  (k_1,\dots,k_m)\in{\cal X}_{k,m}
  \label{defnu}
\end{equation}
where the $\phi_j$ are defined by
\begin{equation}
  \phi_j: k_j \in \N
  \mapsto -\beta k_jv_j + \ln\binom{n_j}{k_j} \in \R\cup\{-\infty\}
\end{equation}
and $Q$ is such that $\nu$ is a probability measure.

The first aim of this paper
is to describe an algorithm
that simulates a sampling according
to $\nu$ in a time that can be bounded
from above by an explicit polynomial
in $k$ and $m$, uniformly over $\beta$,
$(v_j)_{1\leq j \leq m}$ and $(n_j)_{1\leq j\leq m}$.
The reason why we prefer a bound in $k$
and $m$ rather than in the `volume' of the system
$n = \sum_j n_j$, will be clarified later.

A first naive (and wrong) idea to do so
consists in choosing the position (the energy level)
of a first, second,~\dots\ and eventually
$k^{\rm th}$ particle in the following way.
First choose randomly the position
of the first particle according
to the exponential weights associated
with the `free entropies' of 
the empty sites,
that is choose level $j$ with a probability
proportional to $\exp\{-\beta v_j + \ln n_j\}$.
Then decrease by~1 the degeneracy
of the chosen energy level and repeat
the procedure to choose the position
of the second, third,~\dots\ and eventually
$k^{\rm th}$ particle.
It is easy to check that, doing so,
the final distribution of the occupation
numbers $k_1$,~\dots, $k_m$ associated
with the different energy levels,
that is of the numbers of particles
placed in each level,
is in general not given by $\nu$
as soon as $k$ is larger than one.
But it turns out that this naive
idea can be adapted to build an efficient
algorithm to perform approximate samplings
under the Fermi statistics.

Very classically, the fast sampling
performed by the algorithm we will build
will be obtained by running 
a Markov chain $X$ with transition
matrix~$p$ on ${\cal X}_{k,m}$
and with equilibrium measure $\nu$.
The efficiency of the algorithm will be measured
through the bounds that we will be able to give on 
the {\em mixing time}~$t_\epsilon$, defined
for any positive $\epsilon < 1$ by
\begin{eqnarray}
 \label{deft}
 t_\epsilon & := & \inf\{
    t \geq 0 :\:
    d(t) \leq \epsilon
  \}\\
  \label{defd}
  d(t) & := & \max_{\eta\in{\cal X}_{k,m}}\|p^t(\eta,\cdot) - \nu\|_{\rm TV}
\end{eqnarray}
where $\|\cdot\|_{\rm TV}$
stands for the {\em total variation distance}
defined for any probability measures $\nu_1$
and $\nu_2$ on ${\cal X}_{k,m}$ by
\begin{equation}
  \|\nu_1 - \nu_2\|_{\rm TV}:= 
  \max_{A\subset {\cal X}_{k,m}} |\nu_1(A) - \nu_2(A)|
  = \frac{1}{2}\sum_{\eta\in{\cal X}_{k,m}} |\nu_1(\eta) - \nu_2(\eta)|.
\end{equation}
As a consequence,
estimating mixing times is not the only one issue
of this paper, building a `good' Markov chain
is part of the problem.

As far as that part of the problem is concerned,
we propose to build a Metropolis-like algorithm
that uses the `free energies' of the naive approach
to define a conservative dynamics.
Assuming that at time $t\in \N$ the system
is in some configuration $X_t = \eta$
in ${\cal X}_{k,m}$ with $\nu(\eta) > 0$,
and defining for any $\eta = (k_1,\dots,k_m)$
and any distinct $i$ and $j$ in $\{1;\cdots;m\}$
\begin{equation}
  \eta^{ij} := (k'_1,\dots,k'_m)
  \quad\mbox{with}\quad 
  k'_s = \left\{
    \begin{array}{ll}
      k_s & \mbox{for $s \in \{1;\cdots;m\}\setminus\{i;j\}$}\\
      k_i - 1 & \mbox{if $s = i$}\\
      k_j + 1 & \mbox{if $s = j$}
    \end{array}
  \right.,
\end{equation}
the configuration at time $t+1$
will be decided as follows:
\begin{itemize}
\item
  choose a {\em particle} with uniform probability
  (it will stand in a given level $i$ with probability $k_i/k$),
\item 
  choose an {\em energy level} with uniform probability
  (a given level $j$ will be chosen with probability $1/m$),
\item
  with $i$ the level where stood the chosen particle
  and $j$ the chosen energy level, 
  extract a uniform variable $U$ on $[0;1)$
  and set $X_{t+1} = \eta^{ij}$ if $i\neq j$ and
  \begin{equation}
    U < \exp\{
      - \beta v_j + \ln (n_j - k_j) 
      + \beta v_i - \ln (n_i - (k_i - 1))
    \},
  \end{equation}
  $X_{t+1} = \eta$ if not.
\end{itemize}
In other words,
denoting by $[a]^+ = (a + |a|)/2$
the positive part of any real number $a$
and 
with
\begin{equation}
  \psi_j: k_j \in \{0;\cdots;n_j\} 
  \mapsto -\beta v_j + \ln(n_j - k_j) \in \R\cup\{-\infty\},
  \quad j\in\{1;\cdots;m\},
  \label{defpsi1}
\end{equation}
for any distinct $i$ and $j$
in $\{1;\cdots;m\}$
\begin{equation}
  P(X_{t+1} = \eta^{ij} | X_t = \eta) = p(\eta, \eta^{ij})
  = \frac{k_i}{k}\frac{1}{m}\exp\{-[\psi_i(k_i - 1) - \psi_j(k_j)]^+\},
  \label{mela}
\end{equation}
and
\begin{equation}
  P(X_{t+1} = \eta| X_t = \eta)
  = p(\eta, \eta) 
  = 1-\sum_{i\neq j} p(\eta,\eta^{ij}).
\end{equation}
\par\noindent {\bf Remark:}
In order to avoid any ambiguity
in (\ref{mela}) in the case $k_i = 0$,
we set $\psi_i(-1) = +\infty$
(even though the algorithm we described
does not require any convention
for $\psi_i(-1)$).

This Markov chain is certainly irreducible and aperiodic.
To prove that it relaxes to $\nu$
we will check the reversibility of the process
with respect to $\nu$.
Then we will have to
estimate the mixing time of the process.
We will carry out both the tasks in a more 
general setup.

\subsection{A general result}

For any function $f: \N \rightarrow \R$
we define
\begin{equation}
  \nabla^+ f : x \in \N \mapsto \nabla^+_x f := f(x+1) - f(x)
\end{equation}
\begin{equation}
  \nabla^- f : x \in \N\setminus\{0\} \mapsto \nabla^-_x f := f(x-1) - f(x)
\end{equation}
\begin{equation}
  \Delta f : x \in \N\setminus\{0\} \mapsto \Delta_x f := \nabla_x^+f + \nabla_x^-f = -\nabla^-_x(\nabla^+f)
\end{equation}
and we say that a measure $\gamma$ on the integers
\begin{equation}
  \gamma: x\in\N \mapsto e^{\phi(x)}\in\R_+,
\end{equation}
with $\phi: \N\rightarrow\R\cup\{-\infty\}$, is
log-concave if $\N\setminus\gamma^{-1}(\{0\})$
is an interval of the integers and 
\begin{equation}
  \gamma(x)^2\geq\gamma(x-1)\gamma(x+1), \quad x\in\N\setminus\{0\},
\end{equation}
i.e., if $\nabla^+\phi$ is non-increasing, or, equivalently,
$-\Delta\phi$ is non-negative (with the obvious extension
of the previous definitions to such a possibly non-finite $\phi$). 
The measure $\mu$ defined in (\ref{defmu})
is a product of log-concave measures
and the canonical Fermi statistics is such a product
measure normalized over the condition (\ref{canonical}).
\smallskip\par\noindent
{\bf N.B.} {\em From now on, and except for explicit mentioning
  of additional hypotheses, 
  we will only assume that the probability measure $\nu$
  we want to sample
  is a product of log-concave measures normalized
  over the global constraint (\ref{canonical}),
  i.e., that $\nu$ is a probability
  on ${\cal X}_{k,m}$
  that can be written in the form (\ref{defnu})
  with non-increasing $\nabla^+\phi_j$'s.
}
\smallskip\par\noindent
In this more general setup 
we will often refer to the indices $j$ in $\{1;\cdots;m\}$
as {\em sites} rather than energy levels of the system.

Actually the $e^{\phi_j}$'s of the Fermi statistics
are much more than log-concave measures.
They are {\em ultra log-concave} measures
according to the following definition
by Pemantle \cite{Pe} and Liggett \cite{Li}.
\begin{dfn}
  A measure $\gamma:\N\rightarrow\R_+$ is {\em ultra log-concave}
  if $x\mapsto x!\gamma(x)$ is log-concave.
\end{dfn}
In other words $e^{\phi_j}$ is ultra log-concave
if and only if 
\begin{equation}
  \psi_j := \nabla^+\phi_j + \ln(1+\cdot)
  \label{defpsi2}
\end{equation}
is non-increasing (for the Fermi statistics
observe that so are the $\phi_j$'s
and that (\ref{defpsi2}) is consistent
with (\ref{defpsi1})).

For birth and death processes that are reversible
with respect to ultra log-concave measures,
Caputo, Dai Pra and Posta~\cite{CDPP} proved
modified log-Sobolev inequalities
and stronger convex entropy decays,
both giving good upper bounds
on the mixing time of the processes. 
Johnson~\cite{Jo} proved also easier Poincar\'e inequalities
that give weaker bounds on the mixing times.
We refer to~\cite{SC}, \cite{ABCFGMRS}, \cite{MT}
for an introduction to this classical
functional inequality approach to convergence
to equilibrium and we note that for such
birth and death processes the ultra-log concavity
hypotheses allowed for a Bakry-\'Emery like approach
(see~\cite{BE})
to derive (modified) log-Sobolev inequalities.
Actually,~\cite{CDPP} was an attempt to extend
this celebrated analysis to Markov processes with jumps
(see also~\cite{JO} for a more geometric perspective).
But it turns out, as we will soon review,
that beyond the case of birth and death
processes the role of ultra log-concavity
to follow the Bakry-\'Emery line in the discrete setup
is still unclear.
  
In this paper we will not follow the functional analysis
approach to control mixing times.
To bound from above the mixing time of a Markov chain
that is reversible with respect to 
a conditional product
of ultra log-concave measures,
we will follow (and recall in Section~\ref{PT1})
the non less classical, and, in this case, elementary,
probabilistic approach via coalescent coupling.
We will prove: 

\begin{thm}\label{thm1}
  If $\nu$ derives from a product
  of ultra log-concave measures,
  then the Markov chain
  with transition matrix $p$ defined by
  \begin{eqnarray}
  	\label{defp}
    p(\eta,\eta^{ij}) 
    &=& 
    \frac{k_i}{k}\frac{1}{m}\exp\{-[\psi_i(k_i - 1) - \psi_j(k_j)]^+\}, \nonumber\\
    &&
    \qquad \eta = (k_1,\dots,k_m)\in{\cal X}_{k,m}\setminus\nu^{-1}\{0\},\nonumber\\
    &&
    \qquad i\neq j \in \{1;\cdots;m\},\\
    \label{defp2}
    p(\eta, \eta) 
    &=& 1-\sum_{i\neq j} p(\eta,\eta^{ij})
  \end{eqnarray}
  is reversible with respect to $\nu$
  and, for any positive $\epsilon < 1$,
  its mixing time $t_\epsilon$ satisfies
  \begin{equation}
    t_\epsilon \leq km\ln(k/\epsilon).
  \end{equation}
\end{thm}
\noindent
{\bf Proof:} see Section~\ref{PT1}.

The most relevant point of Theorem~\ref{thm1}
with respect to the previous results we know
stands in the uniformity of the upper bound
above the disorder of the system (except
for the ultra log-concavity hypothesis
on $e^{\phi_j}$ in each $j$).
In particular and as far as the Fermi
statistics is concerned, our estimate
does not depend on the temperature,
and, more generally it is independent
from the energy values as well as the
level degeneracies.

To illustrate this fact 
let us start with the case $n_j = 1$
for all $j$.
In this case our dynamics is a simple exclusion
process with site disorder.
Caputo (\cite{Ca1}, \cite{Ca2}) 
proved Poincar\'e inequalities
for such processes, in their continuous time version,
assuming a uniform lower (and upper) bound
on general transition rates, 
while Caputo, Dai Pra and Posta~\cite{CDPP},
looking at particular rates for the process
and still assuming moderate disorder --
that is, uniform lower and upper bounds
on these rates -- proved
a modified log-Sobolev inequality. 
For the particular choice of rates
they made, the upper bound on the mixing time
implied by~\cite{CDPP} 
could not hold in a strong disorder context
(for example with $k=1$, $m = 3$,
$v_1 = v_2 = 0$, $v_3 > 0$ and $\beta\gg 1$).
Our uniformity over the disorder of the system
depends then strongly 
on our particular choice 
for the transition probabilities.
As it is often the case with Markov processes
on discrete state space, the details 
of the dynamics are not less important
than the properties of its equilibrium
measure.

Under the moderate disorder hypotheses of~\cite{Ca1},
\cite{CDPP}, \cite{Ca2}, however,
the upper bound on the mixing time implied by~\cite{CDPP}
and suggested by~\cite{Ca1}, \cite{Ca2} is better than
our upper bound in Theorem~\ref{thm1} (by a factor of order $k$).
But we will see that introducing such moderate disorder
hypotheses in our arguments directly improves our result
by a factor of order $m \geq k$
(see the last remark at the end of Section~\ref{PT1}).

To close the discussion on simple exclusion processes 
we note that those of the arguments in~\cite{Ca2}
that do not depend on the disorder suggest
an upper bound on the mixing time of order $m^2\ln m$.
Theorem~\ref{thm1} improves this estimate when
$k$ is small with respect to $m$.

As far as more general conditional 
non-homogeneous product 
of log-concave measures are concerned,
we stress once again that Theorem~\ref{thm1}
gives a uniform bound over the disorder
that can be improved by a factor of order $m$
by adding moderate disorder hypotheses
(see the last remark
at the end of Section~\ref{last}
and our simulation Appendix).  
Then, such product measures
can be equilibrium measures 
of zero-range processes with a continuous time
generator as in~\cite{BCDPP}, \cite{CDPP}:
\begin{equation}
{\cal L}_\eta f
:=  \frac{1}{m}\sum_{i\neq j}c_i(k_i)\Big(f(\eta^{ij}) - f(\eta)\Big).
\end{equation}
where the $c_j : k_j \in {\mathbb N} \mapsto [0, +\infty)$
are such that $c_j(0) = 0$ and $c_j(k_j) > 0$ for $k_j >0$.  
Indeed, such a process is reversible with respect to $\nu$
provided that 
\begin{equation}
\forall j \in \{1;\cdots;m\},\;
\forall k_j \in {\mathbb N},\;
\phi_j(k_j) = \sum_{0< l \leq k_j} \ln\frac{1}{c_j(l)}\;.
\end{equation}
Boudou, Caputo, Dai Pra and Posta \cite{BCDPP}
proved a Poincar\'e inequality for such a process,
assuming that there exists a positive $c$ such that
\begin{equation}\label{cn1}
\forall j \in \{1;\cdots;m\},\;
\forall k_j \in {\mathbb N},\;
c_j(k_j+1) - c_j(k_j) \geq c.
\end{equation}
This is a moderate disorder hypothesis 
that implies the log-concavity of the $\mu_j$.
To go to modified log-Sobolev inequalities,
i.e., to good mixing time estimates
rather than simple gap estimates,
there is the additional hypothesis in~\cite{CDPP}
that there exists a non-negative $\delta < c$
such that 
\begin{equation}\label{cn2}
\forall j \in \{1;\cdots;m\},\;
\forall k_j \in {\mathbb N},\;
c_j(k_j+1) - c_j(k_j) \leq c + \delta.
\end{equation}
Clearly there are ultra log-concave measures
that do not satisfy~(\ref{cn2}): to have 
an ultra log-concave measure one needs
a strongly decreasing $\nabla^+\phi_j$,
i.e., a strongly increasing $c_j$.
Conversely, it is not true 
that~(\ref{cn1}) and~(\ref{cn2})
imply ultra log-concavity, i.e.,
\begin{equation}\label{ulc}
\forall j \in \{1;\cdots;m\},\;
\forall k_j \geq 1,\;
\ln \frac{k_j+1}{c_j(k_j+1)} \leq \ln \frac{k_j}{c_j(k_j)}\;.
\end{equation}
However, elementary algebra shows
that~(\ref{cn1}) together with~(\ref{cn2}) 
implies 
\begin{equation}
\forall j \in \{1;\cdots;m\},\;
\forall k_j \geq 1,\;
\ln \frac{k_j+1/2}{c_j(k_j+1)} \leq \ln \frac{k_j}{c_j(k_j)}
\end{equation}
and this strangely looks like~(\ref{ulc}).
Therefore we said that the role played
by ultra log-concavity is still unclear
along the functional analysis line of research
in the discrete setup.
We conclude stressing once again
that in this paper we will follow
a different line, that our uniform bound on the mixing time comes from 
an elementary coupling argument and that we do not need
any moderate disorder hypothesis.

 \subsection{Interpolating between sites and particles}

It seems that today available techniques
are such that the less log-concavity we have,
the more homogeneity we need to control 
the convergence to equilibrium.
Staying to the papers mentioned above,
without ultra log-concavity or at least
something that looks like  ultra log-concavity 
we only have Poincar\'e inequalities
for non-homogeneous product of 
log-concave measures,
and without log-concavity
we have modified log-Sobolev inequalities
for homogeneous product measures only
(see \cite{CDPP}).
In addition, the only result we know
for a conservative dynamics
in (weakly) disordered context
and with an equilibrium measure
that is a product of measures  
that are not log-concave
is that of Landim and Noronha Neto~\cite{LNN}
for the (continuous) Ginzburg-Landau process.

We will see that all the ideas of the proof of Theorem~\ref{thm1}
can be extended to deal with a large class
of conditional product of log-concave measures 
that are not ultra log-concave.
To do so, let us define
\begin{equation}
  \delta := 
  \max\left\{
    \lambda \in [0;1] :\:
    \forall j \in \{1;\cdots;m\},\,
    \forall k_j >0,\,
    - \Delta_{k_j}\phi_j \geq \lambda\ln\frac{1 + k_j}{k_j}
  \right\}.
  \label{dfndelta}
\end{equation}
In other words, $\delta$ is the largest
real number in $[0;1]$ for which
all the 
\begin{equation}
  \psi_j^{[\delta]} := \nabla^+\phi_j + \delta\ln(1+\cdot),
  \quad j\in \{1;\cdots;m\},
\end{equation}
are non-increasing.
Denoting by $a\wedge b$ the minimum of two real numbers
$a$, $b$ and defining
\begin{equation}
  l_\delta := k^\delta(k\wedge m)^{1-\delta}
\end{equation}
we will prove:
\begin{thm}\label{thm2}
  If $\delta >0$ 
  then the Markov chain
  with transition matrix $p$ defined by
  \begin{eqnarray}
    p(\eta,\eta^{ij}) 
    &=& 
    \frac{k_i^\delta}{l_\delta}\frac{1}{m}
    \exp\left\{
      -\left[
	  \psi_i^{[\delta]}(k_i - 1) - \psi_j^{[\delta]}(k_j)
      \right]^+
    \right\}, \nonumber\\
    &&
    \qquad \eta = (k_1,\dots,k_m)\in{\cal X}_{k,m}\setminus\nu^{-1}\{0\},\nonumber\\
    &&
    \qquad i\neq j \in \{1;\cdots;m\},\label{defp2a}\\
    p(\eta, \eta) 
    &=& 1-\sum_{i\neq j} p(\eta,\eta^{ij})\label{defp2b}
  \end{eqnarray}
  is reversible with respect to $\nu$
  and, for any positive $\epsilon < 1$,
  its mixing time $t_\epsilon$ satisfies
  \begin{equation}
    t_\epsilon \leq \frac{(k\wedge m)^{1-\delta}}{\delta}km\ln(k/\epsilon).
  \end{equation}
\end{thm}
\noindent
{\bf Remark 1:} By H\"older's inequality, if $\delta>0$
then 
\begin{eqnarray}
  \sum_{i = 1}^{m} k_i^\delta
  = \sum_{i = 1}^{m} k_i^\delta \one^{1-\delta}_{\{k_i \neq 0\}}
  \leq \bigg(\sum_{i = 1}^{m} k_i\bigg)^\delta 
  \bigg(\sum_{i = 1}^{m}\one_{\{k_i\neq 0\}}\bigg)^{1-\delta}
  \leq l_\delta \label{holder}
\end{eqnarray}
and this ensures that~(\ref{defp2a})-(\ref{defp2b})
define a probability matrix.
\par\noindent
{\bf Remark 2:}
As far as this can make sense
in our discrete setup,
we note that the hypothesis
$\delta > 0$ is slightly weaker
than a ``uniform strict log-concavity
hypothesis'' (see~(\ref{dfndelta})).
\par\noindent
{\bf Proof of Theorem 2:} see Section~\ref{PT2}.

For $\delta = 1$ the transition matrix
represents an algorithm starting with a uniform
choice of a particle.
For $\delta = 0$ Theorem~\ref{thm2}
is empty but~(\ref{defp2a})-(\ref{defp2b})
still define a Markov chain $X$ that can
be seen as a particular version
of a discrete state space non-homogeneous
Ginzburg-Landau process. In this case
the transition matrix represents an algorithm
that starts with a uniform choice
of a (non-empty) site.
The case $0 < \delta < 1$ can be seen 
as an interpolation between uniform choices
of site and particle.
More precisely, assuming that at time $t\in \N$
the system is in some configuration
$X_t = \eta = (k_1,\dots,k_m)$
in ${\cal X}_{k,m}\setminus\nu^{-1}\{0\}$,
the configuration at time $t+1$
can be decided as follows.
\begin{itemize}
\item
  Choose a site $i$ or no site at all with
  probabilities proportional to $k_i^\delta$
  and $l_\delta - \sum_i k_i^\delta$. 
\item
  If some site $i$ was chosen, then proceed as
  in the previously described algorithm
  using the functions $\psi_j^{[\delta]}$
  instead of the $\psi_j$'s,
  if not, then set $X_{t+1} = \eta$.
\end{itemize}

\subsection{Last remarks and original motivation}
\label{last}

First, we note that as long as one wants
bounds that are uniform over the disorder,
Theorem 1 gives the right order for
the mixing time: for the Fermi statistic
in the very low temperature regime with,
for example, $v_1 < v_2 < \dots < v_m$
and $n_1, n_2 \geq k$, the equilibrium measure
will be concentrated on $(k,0,\dots,0)$
while, starting from $(0,k,0,\dots,0)$,
the system will reach the ground state
in a time of order $km\ln k$ for large $k$
(this is a coupon-collector estimate).
This fact is illustrated in our simulation
Appendix.

Next, we observe that, with our definitions, $p(\eta, \eta)$
can often be close to one,
especially in strong disorder
situations or when the right and left hand sides
in~(\ref{holder}) are far from each other. 
If one would like to use these results
to perform practical simulations,
then it could be useful to note
that the computational time would still be improved
by implementing an algorithm
that at each step simulates, for a given
configuration $\eta$ on the trajectory
of the Markov chain, the elapsed time 
before the particle reach a {\em different}
configuration (this is a geometric time)
and choose this configuration $\eta'\neq\eta$
according to the (easy to compute) associated law. 
It would then be enough to stop the algorithm
as soon as the total simulated time goes beyond
the mixing time (and then return the last
configuration, that the system occupied
at the mixing time) rather than waiting
for the original algorithm to make a step number
equal to the mixing time.

Turning back to the first naive and wrong
idea, it is interesting to note that it can 
easily be modified to determine the most probable
states of the system, i.e., the configurations 
$\eta^* = (k_1^*,\dots,k_m^*)$ in ${\cal X}_{k,m}$
such that 
\begin{equation}
  \sum_{j = 1}^{m} \phi_j(k_j^*)
  = \max_{k_1+\cdots+k_m = k}  \sum_{j = 1}^{m} \phi_j(k_j).
\end{equation}
One can prove, using the concavity of the $\phi_j$'s,
that the most probable configurations for the system
with $k$ particles can be obtained
from the most probable configurations $(k'_1,\dots, k'_m)$
for the system with $k-1$ particles
simply adding one particle where the corresponding
gain in `free energy' is the highest,
that is in $j^*$ such that
\begin{equation}
  \nabla^+_{k'_{j^*}}\phi_{j^*} = \max_j \nabla^+_{k'_j}\phi_j.
\end{equation}
As a consequence one can place the particles
one by one, each time maximizing this free
energy gain, to build the most probable configuration.

Then, as a referee pointed out, since the sampling problem
is a trivial one in the Poissonian case
(when all the $\mu_j$'s are Poisson measures, so that
$\nu$ is nothing but a multinomial law $M$), one can ask about
the expected time needed to perform a rejection
sampling with respect to the multinomial case.
It is given by $\max_\eta \nu(\eta)/M(\eta)$.
If we take the Fermi statistics with $m = k$
and all the $n_j$ equal to 1, 
then we find (optimizing on the multinomial
parameters by a geometric/arithmetic mean
comparison) $k^k /k! \sim e^k/\sqrt{2\pi k}$.
Of course, the sampling problem for that Fermi
statistics is also a trivial one.
But if we take $\beta = 0$ and all
the $n_j$ equal to 2, it is not anymore a trivial
problem and we find an expected time
for the rejection sampling that is logarithmically
equivalent to $(e/2)^k$.

Turning back to the Fermi statistics we now explain
why we were interested in bounds in $k$ and $m$
rather than in the volume $n$.
Iovanella, Scoppola and Scoppola
defined in~\cite{ISS} an algorithm to individuate
cliques (i.e., complete subgraphs) with $k$
vertices inside a large Erd\"os-Reyni random graph with $n$ vertices.
Their algorithm requires to  
perform repeated approximate samplings
of Fermi statistics in volume $n$, with $k$ particles
and $m = 2k+1$ energy levels.
Now, the key observation is that the largest cliques
in Erd\"os-Reyni graphs with $n$ vertices
are of order $\ln n$, so that $k$ and $m$ in this problem
are logarithmically small with respect to $n$.
Before Theorem~\ref{thm1}
the samplings for their algorithm
were done  by running simple exclusion
processes with $k$ particles 
on the complete graph (with site disorder)
of size $n$.
Such  processes 
converge to equilibrium
in a time of order $kn\ln k$.
Now the samplings are done in a time
of order $km\ln k \sim 2k^2\ln k$,
and that was the original motivation
of the present work.

Finally we note that our bound in $k$ and $m$
is not only useful when $k$ is small with respect
to $n$ but also when $k$ is close to $n$.
In this case we can define a dynamics on the $n-k$
vacancies rather than on the $k$ particles.

\section{Proof of Theorem 1}\label{PT1}

In this section we assume that $\nu$ is a measure on ${\cal X}_{k,m}$
deriving from a product of ultra log-concave measures, which means that we 
can write $\nu(k_1, \ldots, k_m) = (1/Q)\prod_j e^{\phi_j(k_j)}$ and, for 
all $j \in \{1;\cdots;m\}$, $\psi_j$ defined by~(\ref{defpsi2}) is
non-increasing.

\subsection{Reversibility}
\label{reversibility}
We first prove that the transition matrix $p$ defined by~(\ref{defp}) 
and~(\ref{defp2}) is reversible with respect to the measure $\nu$. Let
$\eta = (k_1, \ldots, k_m) \in {\cal X}_{k,m}\setminus\nu^{-1}\{0\}$ 
and let $i\not=j \in \{1;\cdots;m\}$.
We have $p(\eta,\eta^{ij})\neq 0$ if and only if
$\nu(\eta^{ij}) \neq  0$ and, in that case,
\begin{equation}
  \frac{\nu(\eta^{ij})}{\nu(\eta)} = \frac{e^{\phi_i(k_i - 1) + 
  \phi_j(k_j + 1)}}{e^{\phi_i(k_i) + \phi_j(k_j)}}
  = \exp\{\nabla_{k_i}^-\phi_i + \nabla^+_{k_j} \phi_j\}
\end{equation}
while
\begin{eqnarray}
\frac{p(\eta^{ij},\eta)}{p(\eta,\eta^{ij})}
& = & \frac{k_j + 1}{k_i} \exp\{-[\psi_j(k_j) - \psi_i(k_i - 1)]^+ \nonumber\\ 
&&\qquad\qquad\qquad + [\psi_i(k_i - 1) - \psi_j(k_j)]^+\}\\
& = & \exp\{\psi_i(k_i - 1) - \psi_j(k_j) + \ln(k_j+1) - \ln(k_i)\}\\
& = & \exp\{\nabla^+_{k_i-1} \phi_i - \nabla^+_{k_j} \phi_j\}\\
& = &  \exp\{- \nabla_{k_i}^-\phi_i - \nabla^+_{k_j} \phi_j\}
\end{eqnarray}
so that
\begin{equation}
  \nu(\eta)p(\eta,\eta^{ij}) = \nu(\eta^{ij})p(\eta^{ij},\eta).
\end{equation}

\subsection{A few words about the coupling method}

In order to upper bound the mixing time of the Markov chain with transition
matrix $p$, we will use the coupling method. Given a Markov chain $(X^1,X^2)$
on ${\cal X}_{k,m} \times {\cal X}_{k,m}$, we say it is a (Markovian) {\em coupling} 
for the dynamics if both $X^1$ and $X^2$ are Markov chains with transition matrix $p$. 
Given such a coupling, we define the {\em coupling time} $\tau_{\rm couple}$ 
as the first (random) time for which the chains meet, that is
\begin{equation}
  \tau_{\rm couple} := \inf\{t \geq 0 : X^1_t = X^2_t\}.
\end{equation}
In this work, every coupling will also satisfy the condition
\begin{equation}
  \label{ppte}
  t \geq \tau_{\rm couple} \Rightarrow X^1_t = X^2_t.
\end{equation}
Then, it is a well-known fact that for all $t \geq 0$,
\begin{equation}
  \label{idd}
  d(t) \leq \max_{\eta, \theta \in {\cal X}_{k,m}} P(\tau_{\rm couple} > t | X^1_0 = \eta, 
  X^2_0 = \theta)
\end{equation}
where $d(t)$ is defined by~(\ref{defd}).
A proof of this fact as well as an 
exhaustive introduction to mixing time theory can be found in \cite{YP}.

In the proof of both 
\ref{thm1} and \ref{thm2} we will build a coupling
for which there exists a function $\rho$ that measures
in some sense a `distance'
between $X^1_t$ and $X^2_t$
and from which we will get a bound on the mixing time thanks to the 
following proposition.
\begin{prp}
\label{prpdist}
Let $(X^1,X^2)$ be a coupling for a Markov chain with transition matrix $p$. We assume 
that the coupling satisfies the property~(\ref{ppte}). Let $\rho : {\cal X}_{k,m} \times 
{\cal X}_{k,m} \rightarrow \N$ such that $\rho(\eta,\theta) = 0$ if and only if 
$\eta = \theta$.
If $M$ is the maximum of $\rho$ 
and if there exists $\alpha>1$
such that, for all $t \geq 0$,
\begin{equation}
  \label{surmart}
  E\left(\rho(X^1_{t+1}, X^2_{t+1}) | X^1_t, X^2_t\right)
  \leq \left(1-\frac{1}{\alpha}\right) \rho(X^1_t, X^2_t),
\end{equation}
then for all $\epsilon > 0$ the mixing time $t_\epsilon$ of the dynamics is upper
bounded by $\alpha\ln(M/\epsilon)$.
\end{prp}
\par\noindent
{\bf Proof:} Remark that~(\ref{surmart}) actually means that 
$(\rho(X^1_t,X^2_t)(1-1/\alpha)^{-t})_{t\in\N}$ 
is a super-martingale. Taking the expectation we get
\begin{equation}
  E(\rho(X^1_{t+1}, X^2_{t+1})) \leq \left(1-\frac{1}{\alpha}\right) E(\rho(X^1_t, X^2_t))
\end{equation}
As a consequence
\begin{equation}
  \label{majE}
  E(\rho(X^1_t, X^2_t)) \leq \left(1-\frac{1}{\alpha}\right)^t E(\rho(X^1_0, X^2_0)) \leq 
  Me^{-\frac{t}{\alpha}}.
\end{equation}
Since we assume $\rho(\eta,\theta) = 0$ if and only if $\eta = \theta$ and~(\ref{ppte}),
\begin{equation}
  P(\tau_{\rm couple} > t) = P(\rho(X^1_t, X^2_t) > 0) = P(\rho(X^1_t, X^2_t) \geq 1).
\end{equation}
From Markov's inequality and~(\ref{majE}) we deduce 
\begin{equation}
  \label{markov}
  P(\tau_{\rm couple} > t) \leq E(\rho(X^1_t, X^2_t)) \leq Me^{-\frac{t}{\alpha}}.
\end{equation}
Since the upper bound in~(\ref{markov}) is uniform in $X^1_0$ and $X^2_0$, according 
to~(\ref{idd}) we get
\begin{equation}
  d(t) \leq Me^{-\frac{t}{\alpha}}.
\end{equation}
Thus, given $\epsilon > 0$, if $t \geq \alpha\ln(M/\epsilon)$ then $d(t) \leq \epsilon$,
so that $t_\epsilon \leq \alpha\ln(M/\epsilon)$.\qed

\subsection{A colored coupling}

In order to prove Theorem~\ref{thm1}, we introduce a dynamics on the set of the possible
distributions of the $k$ particles in the $m$ energy levels. Therefore we define the set 
of the configurations of $k$ {\em distinguishable} particles in $m$ energy levels
\begin{equation}
  \Omega := \{\omega : \{1;\dots;k\} \to \{1;\dots;m\}\}
\end{equation}
and for every such distribution $\omega$, we define $\xi(\omega) = (k_1, \ldots, k_m) \in 
{\cal X}_{k,m}$ by
\begin{equation}
  \forall i \in \{1;\dots;m\}, k_i := \sum_{x=1}^k \one_{\{\omega(x) = i\}}.
\end{equation}
We will couple two dynamics $\omega^1$ and $\omega^2$ on $\Omega$,
then we will work with the coupling $(\xi(\omega^1),\xi(\omega^2))$.
We will build the coupling $(\omega^1,\omega^2)$
thanks to the coloring we now introduce.

At step $t$, a {\em red-blue coloring} of $(\omega^1_t,\omega^2_t)$ is a couple of functions
$C^1, C^2 : \{1;\dots;k\} \to \{{\rm blue};{\rm red}\}$ such that for all energy level 
$i \in \{1;\dots;m\}$, the number of blue particles in the level $i$ is the same in both
distributions, i.e., if $|A|$ refers to the cardinality of the set $A$,
\begin{equation}
  \label{conditionbleue}
  |\{x : C^1(x) = {\rm blue}, \omega^1_t(x) = i\}| = 
  |\{x : C^2(x) = {\rm blue}, \omega^2_t(x) = i\}|
\end{equation}
and an energy level cannot contain red particles in both distributions
\begin{eqnarray}
  \label{conditionrouge}
  C^1(x) = {\rm red} & \Rightarrow & \forall y, (\omega^2_t(y) = \omega^1_t(x) \Rightarrow
  C^2(y) = {\rm blue}),\\
  C^2(x) = {\rm red} & \Rightarrow & \forall y, (\omega^1_t(y) = \omega^2_t(x) \Rightarrow
  C^1(y) = {\rm blue}).
\end{eqnarray}
As a consequence of~(\ref{conditionbleue}), there exists a one-to-one correspondence 
\begin{equation}
  \label{defPhi}
  \Phi : \{x : C^1(x) = {\rm blue}\} \to \{x : C^2(x) = {\rm blue}\}
\end{equation}
such that for all $x$, $\omega^2_t(\Phi(x)) = \omega^1_t(x)$. 
Since the number 
of blue particles is the same in both distributions,
so is the number of red particles. Moreover, $\xi(\omega^1_t) = \xi(\omega^2_t)$ 
if and only if all particles are blue. Then we are willing to provide a coupling 
$(\omega^1,\omega^2)$ for which the number of red particles is non-increasing.
This number does not depend on the red-blue coloring
and at any time $t$ we will call it $\rho_t$. 
Writing 
$\xi(\omega^1_t) = (k^1_1, \ldots, k^1_m)$ and $\xi(\omega^2_t) = (k^2_1, \ldots, k^2_m)$
we have the identities
\begin{equation}
\rho_t = \frac{1}{2}\sum_{i=1}^m|k^1_i - k^2_i| = \sum_{i=1}^m[k^1_i-k^2_i]^+ 
= \sum_{i=1}^m[k^2_i-k^1_i]^+.
\end{equation}

We are now ready
to build our coupling $(\omega^1,\omega^2)$.
Given $(\omega^1_t,\omega^2_t)$ and a red-blue coloring $(C^1,C^2)$ 
at step $t$ (such a coloring 
certainly exists for any couple of distributions $(\omega^1_t,\omega^2_t)$), let $x^1$ 
be a uniform random integer in $\{1;\dots;k\}$.
\begin{itemize}
  \item If $C^1(x^1) = {\rm blue}$: then we set $x^2 = \Phi(x^1)$ where $\Phi$ is 
  provided by~(\ref{defPhi}).
  \item If $C^1(x^1) = {\rm red}$: let $x^2$ be a uniform random integer in $\{x :
  C^2(x) = {\rm red}\}$. 
\end{itemize}
\begin{lmm}
\label{lmmx2}
The random integer $x^2$ has a uniform distribution on the set $\{1;\dots;k\}$.
\end{lmm}
\par\noindent
{\bf Proof:} We write 
\begin{eqnarray}
  P(x^2 = y) & = & \sum_{x=1}^{k} P(x^2 = y | x^1 = x)P(x^1=x) \\
  & = & \sum_{x : C^1(x) = {\rm blue}} \one_{\{y = \Phi(x)\}}\frac{1}{k} \nonumber\\
  && +\sum_{x : C^1(x) = {\rm red}} \frac{\one_{\{C^2(y) = {\rm red}\}}}{\rho_t} \frac{1}{k}\\
  & = & \frac{\one_{\{C^2(y) = {\rm blue}\}}}{k} +
        \frac{\one_{\{C^2(y) = {\rm red}\}}}{\rho_t}\frac{\rho_t}{k} \\
  & = & \frac{1}{k}.
\end{eqnarray}
\qed

\medskip\par\noindent
We then choose an energy level $j$ uniformly in $\{1;\dots;m\}$ and we write 
$\xi(\omega^1_t) = (k^1_1, \ldots, k^1_m)$ and $\xi(\omega^2_t) = (k^2_1, \ldots, k^2_m)$.

\begin{itemize}
\item 
  If $C^1(x^1) = {\rm blue}$: 
  then $C^2(x^2) = {\rm blue}$, and $x^1$, $x^2$ are
  in the same energy level $i := \omega^1_t(x^1) = \omega^2_t(x^2)$. Let $a\not=b \in
  \{1;2\}$ such that $\psi_i(k^a_i - 1) - \psi_j(k^a_j) \leq \psi_i(k^b_i - 1) - 
  \psi_j(k^b_j)$. Then, 
  \begin{eqnarray}
    p_a & := & \exp\{-[\psi_i(k^a_i-1) - \psi_j(k^a_j)]^+\}\\
    & \geq &
    \exp\{-[\psi_i(k^b_i-1) - \psi_j(k^b_j)]^+\} \;=:\; p_b.
  \end{eqnarray}
  Let $U$ be a uniform random variable on $[0;1)$.
    \begin{itemize}
    \item 
      If $U < p_b$:
      then we set $\omega^a_{t+1}(x) = \omega^a_t(x)$ for all $x \not=
      x^a$, $\omega^a_{t+1}(x^a) = j$, $\omega^b_{t+1}(x) = \omega^b_t(x)$ for all $x \not=
      x^b$ and $\omega^b_{t+1}(x^b) = j$.
      Then $\xi(\omega^a_{t+1}) = (\xi(\omega^a_t))^{ij}$ 
      and $\xi(\omega^b_{t+1}) = (\xi(\omega^b_t))^{ij}$, and for any red-blue coloring of
      $(\omega^1_{t+1},\omega^2_{t+1})$
      the number of red particles $\rho_{t+1}$ remains the same
      (since both particles $x^1$ and $x^2$ have moved together).
    \item 
      If $p_b \leq U < p_a$: then we set $\omega^a_{t+1}(x) = \omega^a_t(x)$ for all 
      $x \not= x^a$, $\omega^a_{t+1}(x^a) = j$
      and $\omega^b_{t+1}(x) = \omega^b_t(x)$ for all $x$.
      Then $\xi(\omega^a_{t+1}) = (\xi(\omega^a_t))^{ij}$ and $\xi(\omega^b_{t+1}) = 
      \xi(\omega^b_t)$.
      The only situation in which the number of red particles could increase
      is the following: $k^a_j \geq k^b_j$ and $k^a_i \leq k^b_i$.
      Since $\psi_i$ and $\psi_j$
      are non-increasing, this would imply 
      $p_a \leq p_b$ that contradicts $p_b \leq U < p_a$. Then, $\rho_{t+1} \leq \rho_t$ 
      for any red-blue coloring of $(\omega^1_{t+1},\omega^2_{t+1})$.
    \item 
      If $U \geq p_a$: then we set $\omega^a_{t+1}=\omega^a_t$ and 
      $\omega^b_{t+1}=\omega^b_t$. Then $\xi(\omega^a_{t+1}) = \xi(\omega^a_t)$, 
      $\xi(\omega^b_{t+1}) = \xi(\omega^b_t)$ and $\rho_{t+1}=\rho_t$.
    \end{itemize}
\item 
  If $C^1(x^1) = {\rm red}$: then $C^2(x^2) = {\rm red}$, and we define $i^1 := 
  \omega^1_t(x^1)$ and $i^2 := \omega^2_t(x^2)$. Note that according 
  to~(\ref{conditionrouge}), $i^1 \not= i^2$.
  Let then $V^1, V^2$ be two independent uniform
  random variables on $[0;1)$. 
    \begin{itemize}
    \item 
      If $V^1 < \exp\{-[\psi_{i^1}(k^1_{i^1}-1) - \psi_j(k^1_j)]^+\}$ then we set 
      $\omega^1_{t+1}(x) = \omega^1_t(x)$ for all $x \not= x^1$
      and $\omega^1_{t+1}(x^1) = j$
      and then $\xi(\omega^1_{t+1}) = (\xi(\omega^1_t))^{ij}$; 
      otherwise we leave $\omega^1_{t+1} = \omega^1_t$.
    \item 
      If $V^2 < \exp\{-[\psi_{i^2}(k^2_{i^2}-1) - \psi_j(k^2_j)]^+\}$ then we set 
      $\omega^2_{t+1}(x) = \omega^2_t(x)$ for all $x \not= x^2$
      and $\omega^2_{t+1}(x^2) = j$
      and then $\xi(\omega^2_{t+1}) = (\xi(\omega^2_t))^{ij}$; 
      otherwise we leave $\omega^2_{t+1} = \omega^2_t$.
    \end{itemize}
    Whether red particles move or not,
    the number of red particles cannot increase, so it is clear
    that $\rho_{t+1} \leq \rho_t$.
\end{itemize}
We conclude
\begin{prp}
  $(\rho_t)_{t\in\N}$ is a non-increasing process.
\end{prp}
and claim
\begin{prp}
The process $(\xi(\omega^1),\xi(\omega^2))$
is a coupling for the Markov chain with transition matrix $p$.
\end{prp}
\par\noindent
{\bf Proof:} Writing $\xi(\omega^1_t) = \eta$ and given $i,j \in \{1;\dots;m\}$,
\begin{eqnarray}
  && P(\xi(\omega^1_{t+1}) = \eta^{ij}) \nonumber\\
  && \quad = \quad P(\xi(\omega^1_{t+1}) = \eta^{ij} | C^1(x^1) = {\rm blue})
  P(C^1(x^1) = {\rm blue})\nonumber\\
  && \quad\quad\quad +\; P(\xi(\omega^1_{t+1}) = \eta^{ij} | C^1(x^1) = {\rm red})
  P(C^1(x^1) = {\rm red}).
\end{eqnarray}
We have
\begin{equation}
  P(C^1(x^1) = {\rm red}) = 1 - P(C^1(x^1) = {\rm blue}) = \frac{\rho_t}{k}
\end{equation}
and
\begin{eqnarray}
  && P(\xi(\omega^1_{t+1}) = \eta^{ij} | C^1(x^1) = {\rm blue}) \nonumber\\
  &&\quad = \quad
  P(\omega^1_t(x^1) = i | C^1(x^1) = {\rm blue})
  \times \frac{1}{m} \times P(U<p_1)\\
  &&\quad = \quad
  \frac{k^1_i \wedge k^2_i}{(k-\rho_t)m}\exp\{-[\psi_i(k^1_i-1) - \psi_j(k^1_j)]^+\}\\
  \nonumber\\
  &&  P(\xi(\omega^1_{t+1}) = \eta^{ij} | C^1(x^1) = {\rm red}) \nonumber\\
  &&\quad = \quad
  P(\omega^1_t(x^1) = i | C^1(x^1) = {\rm red})
  \times \frac{1}{m} \times P(V^1<p_1)\\
  &&\quad = \quad
  \frac{[k^1_i-k^2_i]^+}{\rho_tm}\exp\{-[\psi_i(k^1_i-1) - \psi_j(k^1_j)]^+\}
\end{eqnarray}
which finally leads to 
\begin{equation}
  \label{xi}
  P(\xi(\omega^1_{t+1}) = \eta^{ij})
  = \frac{k^1_i}{km}\exp\{-[\psi_i(k^1_i-1) - \psi_j(k^1_j)]^+\}
  = p(\eta,\eta^{ij}).
\end{equation}
Then, $P(\xi(\omega^1_{t+1}) = \eta^{ij} | \omega^1_t, \omega^2_t)$ 
depends on $\xi(\omega^1_t)$ only,
which means that $\xi(\omega^1)$ is a Markov chain. 
Besides, according to~(\ref{xi}) its transition matrix is $p$. 
Finally, by Lemma~\ref{lmmx2}, the same is true
for $\xi(\omega_2)$.
\qed

\medskip\par\noindent
From now on, we will write $X^1_t := \xi(\omega^1_t)$ and $X^2_t := \xi(\omega^2_t)$.
The previous proposition ensures that $(X^1,X^2)$ 
is a coupling for the dynamics with transition 
matrix $p$.

\subsection{Estimating the coupling time}

We will use Proposition~\ref{prpdist}.
Since $(\rho_t)_{t\in\N}$ is non increasing
and $\rho_t = 0$ if and only if $X^1_t = X^2_t$,
all we have to do is to estimate from below
the probability of the event $\{\rho_{t+1} < \rho_t\}$.

\begin{prp}
\label{bonchoix}
If at step $t$ of the coupled dynamics $(\omega^1,\omega^2)$ we assume that red particles
have been chosen, i.e., $C^1(x^1) = {\rm red} = C^2(x^2)$, then there is a choice of 
$j \in \{i^1;i^2\}$ for which the number of red particles decreases with probability 1, where
$i^1$ (resp. $i^2$) still refers to $\omega^1_t(x^1)$ (resp. $\omega^2_t(x^2)$).
\end{prp}
\par\noindent
{\bf Proof:} If the inequalities
\begin{eqnarray}
\label{ineq1}\exp\{-[\psi_{i^1}(k^1_{i^1}-1) - \psi_{i^2}(k^1_{i^2})]^+\} & < & 1\\
\label{ineq2}\exp\{-[\psi_{i^2}(k^2_{i^2}-1) - \psi_{i^1}(k^2_{i^1})]^+\} & < & 1
\end{eqnarray}
hold together, then $\psi_{i^1}(k^1_{i^1}-1) > \psi_{i^2}(k^1_{i^2})$ 
and $\psi_{i^2}(k^2_{i^2}-1) > \psi_{i^1}(k^2_{i^1})$. Besides,
since $C^1(x^1) = {\rm red}$, according to~(\ref{conditionrouge}) $k^1_{i^1} > k^2_{i^1}$
from which we get $k^1_{i^1} - 1 \geq k^2_{i^1}$ and,
since $\psi_{i^1}$ is non-increasing,
$\psi_{i^1}(k^1_{i^1} - 1) \leq \psi_{i^1}(k^2_{i^1})$.
Likewise, since $C^2(x^2) = {\rm red}$
we have $\psi_{i^2}(k^2_{i^2} - 1) \leq \psi_{i^2}(k^1_{i^2})$. We finally may write
\begin{equation}
\psi_{i^1}(k^1_{i^1}-1) 
> \psi_{i^2}(k^1_{i^2})
\geq \psi_{i^2}(k^2_{i^2} - 1)
> \psi_{i^1}(k^2_{i^1})
\geq \psi_{i^1}(k^1_{i^1} - 1)
\end{equation}
which is absurd. As a result, either~(\ref{ineq1}) or~(\ref{ineq2}) is false.
For instance let
us assume that~(\ref{ineq1}) is false,
then if $j=i^2$, with probability 1 we have
$\omega^1_{t+1}(x) = \omega^1_t(x)$
for all $x\not=x^1$, $\omega^1_{t+1}(x^1) = i^2 \not= \omega^1_t(x^1)$
and $\omega^2_{t+1} = \omega^2_t$.
Then, the number of red particles for any red-blue coloring of
$(\omega^1_{t+1},\omega^2_{t+1})$ is exactly $\rho_{t+1} = \rho_t - 1$. 
\qed

\medskip\par\noindent
It follows
\begin{cor} 
\label{zucchero}
At step $t$, if $\rho_t > 0$,
the probability for $\rho_{t+1}$ to be $\rho_t - 1$ is at least $\rho_t/km$.
\end{cor}

Consequently, and owing to the fact that $\rho_t$ cannot increase,
we have the inequality 
$E(\rho_{t+1}-\rho_t|X^1_t,X^2_t) \leq -P(\rho_{t+1}=\rho_t-1|X^1_t,X^2_t)$ 
from which we deduce
\begin{equation}
E(\rho_{t+1}|X^1_t,X^2_t) \leq \left(1-\frac{1}{km}\right)\rho_t.
\end{equation}
Then we can apply Proposition~\ref{prpdist}
to $\rho(X^1_t,X^2_t)=\rho_t$ with $M=k$ and $\alpha=km$,
which finally proves Theorem~\ref{thm1}.

\medskip\par\noindent
{\bf Remark:} Adding ``moderate disorder hypotheses''
we can gain a lot, depending on the specific model we consider.
For example in the case of our simple exclusion dynamics,
i.e., for the Markov chain associated with the Fermi statistics
when all the $n_j$'s are equal to 1, we can gain a factor of order $m$.
We can indeed assume that $k/m \leq 1/2$ (if not, we consider a dynamics
on the vacancies rather than on the particles
as mentioned in Section~\ref{last}) and, if, for example, $\nu$
derives from an homogeneous product measure, i.e., $\psi_i = \psi_j$
for all $i$ and $j$, then there are much more than one choice
for $j$ for which the number of red particles will decrease with
probability 1 when red particles were chosen: any of the more than
$m/2$ vacant sites in $X^1_t$ or $X^2_t$ will do the job.
In this case we gain a factor $m/2$, first in Corollary~\ref{zucchero},
then in the final estimate.
More generally, if the functions $|\psi_i - \psi_j|$ are uniformly
bounded (this is a moderate disorder hypothesis), then, when red particles
are chosen together with one of these more than $m/2$ vacant sites,
the probability that the number of red particles decrease is bounded
away from zero: once again we gain a factor of order $m$.

\section{Proof of Theorem 2}\label{PT2}

We now work with the dynamics defined by (\ref{defp2a})-(\ref{defp2b})
and we assume $\delta > 0$.
First, the reversibility of this dynamics with respect to $\nu$ still holds,
with exactly the same 
computation as in 
Subsection~\ref{reversibility}.
However, it is no longer possible to work
with an underlying process $\omega_t \in \Omega$ since the factor $k_i^{\delta}$ cannot
stand for a number of particles as soon as $\delta < 1$. Therefore we need to adapt the 
coupling $(X^1,X^2)$ directly on ${\cal X}_{k,m}$.

\subsection{Generalizing the previous coupling}

At step $t$, let us assume $X^1_t = (k^1_1, \ldots, k^1_m) \in {\cal X}_{k,m}$ and
$X^2_t = (k^2_1, \ldots, k^2_m) \in {\cal X}_{k,m}$. We define the following sets:
\begin{eqnarray}
R_1 & := & \{i \in \{1;\dots;m\} : k^1_i > k^2_i\},\\
R_2 & := & \{i \in \{1;\dots;m\} : k^2_i > k^1_i\},\\
B   & := & \{i \in \{1;\dots;m\} : k^1_i = k^2_i\},\label{defB}
\end{eqnarray}
and the following quantities:
\begin{eqnarray}
    w_1 :=  \sum_{i \in R_1} (k^1_i)^{\delta}, 
    &&
    w_1':=  \sum_{i \in R_2} (k^1_i)^{\delta},\\
    w_2 :=  \sum_{i \in R_2} (k^2_i)^{\delta},
    &&
    w_2':=  \sum_{i \in R_1} (k^2_i)^{\delta},
\end{eqnarray}
\begin{equation}
    w_B := \sum_{i \in B} (k^1_i)^{\delta}
    = \sum_{i \in B} (k^2_i)^{\delta}.
\end{equation}
Finally we define
\begin{equation}
  \rho_t := \frac{1}{2}\sum_{i=1}^m|k^1_i - k^2_i| = \sum_{i=1}^m[k^1_i-k^2_i]^+ 
  = \sum_{i=1}^m[k^2_i-k^1_i]^+.
\end{equation}

Keeping in mind the previous coloring, $R_1$ (resp. $R_2$) is the set of sites in 
which there are red particles for the first (resp. the second) configuration, $B$ is 
the set of sites in which there are only blue particles 
or no particles for both configurations, $w_1$
(resp. $w_2$) is proportional to the probability to choose a site for the first (resp.
the second) configuration in which there are red particles, $w_1'$ (resp. $w_2'$) is 
proportional to the probability to choose a site for the first (resp. the second) 
configuration in which there are only blue particles
while there are red particles in the second 
(resp. the first) configuration, $w_B$ is proportional to the 
probability to choose a site in which there are only blue particles for both 
configurations, and $l_{\delta} - (w_B + w_1 + w_1')$ (resp. $l_{\delta} - (w_B + 
w_2 + w_2')$) is proportional to the probability not to choose any site for the first
(resp. the second) configuration, and it can be positive as soon as $\delta < 1$.
Finally, $\rho_t$ still stands for the number of
red particles and it is clear that $\rho_t = 0$ if and only if $X^1_t = X^2_t$.
\smallskip\par\noindent
{\bf N.B.} {\em In the remaining part of this subsection, we assume $w_1+w_1'\geq w_2+w_2'$
in order to not overload the notations and not increase the number of cases to 
investigate. Obviously, the case $w_1+w_1' \leq w_2+w_2'$ is exactly symmetric.}
\smallskip\par\noindent
Let $I$ be a uniform random variable on $[0;l_{\delta})$.
\begin{itemize}
  \item[(i)] If $I < w_B$: then there exists a unique $i \in B$ such that
    \begin{equation}
      \sum_{i' \in B ; i' < i}(k^1_{i'})^{\delta} \leq I 
      < \sum_{i' \in B ; i' \leq i}(k^1_{i'})^{\delta}
    \end{equation}
    and we set $i^1 = i^2 = i$.
    \par\noindent\underline{Remark:} We then have, for all $i \in \{1;\dots;m\}$,
    \begin{eqnarray}
      P(i^1=i|{\rm (i)}) & = & \one_{\{i \in B\}}\frac{(k^1_i)^{\delta}}{w_B},\\
      P(i^2=i|{\rm (i)}) & = & \one_{\{i \in B\}}\frac{(k^1_i)^{\delta}}{w_B} = 
      \one_{\{i \in B\}}\frac{(k^2_i)^{\delta}}{w_B}.
    \end{eqnarray}
  \item[(ii)] If $w_B \leq I < w_B + w_1'$: then there exists a unique $i \in R_2$ such that
    \begin{equation}
      \sum_{i' \in R_2 ; i' < i}(k^1_{i'})^{\delta} \leq I - w_B 
      < \sum_{i' \in R_2 ; i' \leq i}(k^1_{i'})^{\delta}
    \end{equation}
    and we set $i^1 = i^2 = i$.
    \par\noindent\underline{Remark:} We then have, for all $i \in \{1;\dots;m\}$,
    \begin{eqnarray}
      P(i^1=i|{\rm (ii)}) & = & \one_{\{i \in R_2\}}\frac{(k^1_i)^{\delta}}{w_1'},\\
      P(i^2=i|{\rm (ii)}) & = & \one_{\{i \in R_2\}}\frac{(k^1_i)^{\delta}}{w_1'}.
    \end{eqnarray}
  \item[(iii)] If $w_B + w_1' \leq I < w_B + w_1' + w_1$:
    then there exists a unique $i \in R_1$
    such that 
    \begin{equation}
      \sum_{i' \in R_1 ; i' < i}(k^1_{i'})^{\delta} \leq I - w_B - w_1' 
      < \sum_{i' \in R_1 ; i' \leq i}(k^1_{i'})^{\delta}.
    \end{equation}
    We set $i^1 = i$ and we define $u := I - w_B - w_1' - \sum_{i' \in R_1 ; i' < i}
    (k^1_{i'})^{\delta}$, so that $0 \leq u < (k^1_i)^{\delta}$. Notice that since 
    $i \in R_1$, $k^1_i > k^2_i$. If $u < (k^2_i)^{\delta}$ then we set $i^2 = i$. 
    Otherwise, for all $i' \in R_2$ we write $\overline{k}_{i'} := \sum_{\ell \in R_2; \ell<i'}
    (k^1_{\ell})^{\delta}$ and we denote by $T$ the disjoint union of intervals
    \begin{equation}
      T := \left\{\bigcup_{i' \in R_2}\left[\overline{k}_{i'} + (k^1_{i'})^{\delta} ; \overline{k}_{i'} + 
	(k^2_{i'})^{\delta}\right)\right\} \cup [w_2 + w_2' ; w_1 + w_1').
    \end{equation}
    Let $I'$ be a uniform random variable on $T$.
    If there exists $i' \in R_2$ such that
    $I'\in \left[\overline{k}_{i'} + (k^1_{i'})^{\delta} ; \overline{k}_{i'} + (k^2_{i'})^{\delta}\right)$
    then
    we set $i^2 = i'$. Else we do not define $i^2$.
    \par\noindent\underline{Remark:} We then have, for all $i \in \{1;\dots;m\}$,
    \begin{equation}
      P(i^1=i|{\rm (iii)}) = \one_{\{i \in R_1\}}\frac{(k^1_i)^{\delta}}{w_1}
    \end{equation}
    and for all $i' \in \{1;\dots;m\}$,
    \begin{eqnarray}
      &&P(i^2=i'|{\rm (iii)})\nonumber\\
      && \quad = \quad 
      \one_{\{i' \in R_1\}}\frac{(k^2_{i'})^{\delta}}{w_1}\nonumber\\
      && \quad\quad\quad
      +\; \one_{\{i' \in R_2\}}\sum_{i \in R_1}\frac{(k^1_i)^{\delta}-(k^2_i)^{\delta}}{w_1}
      \times \frac{(k^2_{i'})^{\delta}-(k^1_{i'})^{\delta}}{\lambda(T)}\label{iii2}
    \end{eqnarray}
    where 
    \begin{eqnarray}
      \lambda(T) & := & \sum_{i' \in R_2}\left((k^2_{i'})^{\delta}-(k^1_{i'})^{\delta}\right)
      + (w_1 + w_1') - (w_2 + w_2')\\
      & = & w_1 - w_2'\label{lambdaT}
    \end{eqnarray}
    so that
    \begin{equation}
      P(i^2=i'|{\rm (iii)}) = \one_{\{i' \in R_1\}}\frac{(k^2_{i'})^{\delta}}{w_1} +
      \one_{\{i' \in R_2\}}\frac{(k^2_{i'})^{\delta}-(k^1_{i'})^{\delta}}{w_1}.
    \end{equation}
  \item[(iv)] 
    If $I \geq w_B + w_1 + w_1'$ (this case 
    cannot occur when $\delta = 1$):
    then we do not define $i^1$ and $i^2$.
\end{itemize}
Before going ahead with the definition of our coupling
we note, as a direct consequence of our remarks in (i), (ii), (iii)
and of of the fact that $I$ has a uniform distribution:
\begin{prp}
  \label{choixi}
  For all $i \in \{1;\dots;m\}$, $P(i^1 = i | X^1_t, X^2_t) = (k^1_i)^{\delta}/l_{\delta}$ and
  $P(i^2 = i | X^1_t, X^2_t) = (k^2_i)^{\delta}/l_{\delta}$.
\end{prp}

\medskip\par\noindent
We then choose an integer $j \in \{1;\dots;m\}$ with uniform law
and we distinguish once again between our four previous cases.
\begin{itemize}
\item[(i)] 
  If $i^1 \in B$: then $i^2=i^1$, we just write $i^1=i^2=i$. Then $k^1_i = k^2_i$. 
  Thus, let $a,b \in \{1;2\}$ such that $a\not=b$ and $k^a_j \leq k^b_j$. Since both 
  $\psi_i^{[\delta]}$ and $\psi_j^{[\delta]}$ are non-increasing, 
  \begin{eqnarray}
    p_a & := & \exp\left\{-\left[
      \psi_i^{[\delta]}(k^a_i-1) - \psi_j^{[\delta]}(k^a_j)
   \right]^+\right\}\\
    & \geq &
    \exp\left\{-\left[
      \psi_i^{[\delta]}(k^b_i-1) - \psi_j^{[\delta]}(k^b_j)
    \right]^+\right\} \;=:\; p_b.
  \end{eqnarray}
  Let $U$ be a uniform random variable on $[0;1)$.
  \begin{itemize}
  \item 
    If $U < p_b$: we set $X^a_{t+1} = (X^a_t)^{ij}$ and $X^b_{t+1} = (X^b_t)^{ij}$.
  \item 
    If $p_b \leq U < p_a$: we set $X^a_{t+1} = (X^a_t)^{ij}$ and $X^b_{t+1} = X^b_t$.
  \item 
    If $p_a \leq U$: we set $X^a_{t+1} = X^a_t$ and $X^b_{t+1} = X^b_t$.
  \end{itemize}
  In any of these cases, we certainly have $\rho_{t+1} = \rho_t$.
\item[(ii)] 
  If $i^1 \in R_2$: then $i^2=i^1$, we just write $i^1=i^2=i$.
  Let $a,b \in \{1;2\}$ such that $a\not=b$ 
  and $\psi_i^{[\delta]}(k^a_i-1) - \psi_j^{[\delta]}(k^a_j) 
  \leq
  \psi_i^{[\delta]}(k^b_i-1) - \psi_j^{[\delta]}(k^b_j)$,
  so that
  \begin{eqnarray}
    p_a & := & \exp\left\{-\left[
      \psi_i^{[\delta]}(k^a_i-1) - \psi_j^{[\delta]}(k^a_j)
    \right]^+\right\}\\
    & \geq &
    \exp\left\{-\left[
      \psi_i^{[\delta]}(k^b_i-1) - \psi_j^{[\delta]}(k^b_j)
    \right]^+\right\} \;=:\; p_b.
  \end{eqnarray}
  Let $U$ be a uniform random variable on $[0;1)$.
  \begin{itemize}
  \item 
    If $U < p_b$: we set $X^a_{t+1} = (X^a_t)^{ij}$ and $X^b_{t+1} = (X^b_t)^{ij}$.
  \item 
    If $p_b \leq U < p_a$: we set $X^a_{t+1} = (X^a_t)^{ij}$ and $X^b_{t+1} = X^b_t$.
  \item 
    If $p_a \leq U$: we set $X^a_{t+1} = X^a_t$ and $X^b_{t+1} = X^b_t$.
  \end{itemize}
  In the last case we obviously have $\rho_{t+1} = \rho_t$.
  In the first case the particles move together and $\rho_{t+1}=\rho_t$.
  In the second case the number of red particles could increase
  only if $k_i^a\leq k_i^b$ and $k_j^a \geq k_j^b$,
  but, since $\psi_i^{[\delta]}$ and $\psi_j^{[\delta]}$
  are non increasing, this would contradict $p_a > p_b$.
  As a consequence we have $\rho_{t+1} \leq \rho_t$
  in all the three cases. 
\item[(iii)] 
  If $i^1 \in R_1$: there are three cases for $i^2$.
  Either $i^2=i^1=i$ and this case
  is the symmetric of $(ii)$.
  Or $i^2$ is randomly chosen in $R_2$, and we define
  \begin{eqnarray}
    p_1 := \exp\left\{-\left[
      \psi_{i^1}^{[\delta]}(k^1_{i^1}-1) - \psi_j^{[\delta]}(k^1_j)
    \right]^+\right\},
    \label{dfnp1}\\
    p_2 := \exp\left\{-\left[
      \psi_{i^2}^{[\delta]}(k^2_{i^2}-1) - \psi_j^{[\delta]}(k^2_j)
    \right]^+\right\}.   
  \end{eqnarray}
  Or else $i^2$ is not defined,
  and we set $p_2 := 0$ still defining $p_1$ by~(\ref{dfnp1}).
  Let then $V^1$, $V^2$ be independent uniform random variables on $[0;1)$.
  \begin{itemize}
  \item 
    If $V^1 < p_1$: we set $X^1_{t+1} = (X^1_t)^{i^1j}$, else we set $X^1_{t+1} = X^1_t$.
  \item 
    If $V^2 < p_2$: we set $X^2_{t+1} = (X^2_t)^{i^2j}$, else we set $X^2_{t+1} = X^2_t$.
  \end{itemize}
  In the first case we have $\rho_{t+1}\leq \rho_t$ as previously.
  In the last two cases we also have $\rho_{t+1}\leq \rho_t$ since only particles
  from $R_1$ in the first configuration and from $R_2$ in the second configuration
  can move.
\item[(iv)]
  If $i_1$ and $i_2$ are not defined: then we simply set
  $(X^1_{t+1}, X^2_{t+1}) = (X^1_{t}, X^2_{t})$
  and we have $\rho_{t+1} = \rho_t$. 
\end{itemize}

In any of the previous cases,
once $i^1$, $i^2$ and $j$ have been defined,
the probability for
$X^1_{t+1}$ (resp. $X^2_{t+1}$) to be $(X^1_t)^{i^1j}$ (resp. $(X^2_t)^{i^2j}$) is 
$\exp\{-[\psi_{i^1}^{[\delta]}(k^1_{i^1}-1) - \psi_j^{[\delta]}(k^1_j)]^+\}$
(resp. $\exp\{-[\psi_{i^2}^{[\delta]}(k^2_{i^2}-1) - \psi_j^{[\delta]}(k^2_j)]^+\}$).
Thus, according to Proposition~\ref{choixi},
the fact that $j$ is uniformly chosen in $\{1;\dots;m\}$
and our study on the variation of $\rho$
we conclude:
\begin{prp}
The process $(X^1,X^2)$ is a coupling for the dynamics defined by 
(\ref{defp2a})-(\ref{defp2b}) and such that
$(\rho_t)_{t\in\N}$ is non-increasing.
\end{prp}

\subsection{Estimating the coupling time}

Similarly to the proof of Theorem~\ref{thm1} we will use
Proposition~\ref{prpdist}:
since $(\rho_t)_{t\in\N}$ is non-increasing
it will be enough to give a lower bound 
for the probability of $\{\rho_{t+1} < \rho_t\}$. 

\begin{prp}
\label{choixj2}
If at step $t$ of the coupled dynamics $(X^1,X^2)$,
we assume that ``red particles have been chosen'',
i.e., $i^1 \in R_1$ and $i^2 \in R_2$,
then, there is a choice of $j \in \{i^1;i^2\}$ for which 
the number of red particles decreases with probability 1.
\end{prp}
\par\noindent
{\bf Proof:} Assuming $i^1 \in R_1$ and $i^2 \in R_2$ yields 
$k^1_{i^1} > k^2_{i^1}$ and $k^1_{i^2}< k^2_{i^2}$.
Using exactly the same argument as for Proposition~\ref{bonchoix} we prove that either
$\exp\{-[\psi_{i^1}^{[\delta]}(k^1_{i^1}-1) - \psi_{i^2}^{[\delta]}(k^1_{i^2})]^+\} = 1$ 
or 
$\exp\{-[\psi_{i^2}^{[\delta]}(k^2_{i^2}-1) - \psi_{i^1}^{[\delta]}(k^2_{i^1})]^+\} = 1$.
Eventually,
if one red particle in some configuration  
moves to a site with a red particle in the other configuration,
then both particles turn blue
and the number of red particles decreases by one.
\qed

\begin{cor}
At step $t$, the probability for $\rho_{t+1}$ to be $\rho_t - 1$ is at least $\delta 
k^{\delta-1}\rho_t/ml_{\delta}$.
\end{cor}
\par\noindent
{\bf Proof:} The probability to choose $i^1 \in R_1$ and $i^2 \in R_2$ is 
\begin{eqnarray}
  &&P(i^1 \in R_1, i^2 \in R_2)\nonumber\\ 
  && \quad = \quad \sum_{i'\in R_2}\sum_{i\in R_1} P(i^2 = i'|i^1=i)\times P(i^1=i)\\
  && \quad = \quad \sum_{i'\in R_2}\sum_{i\in R_1} 
  \frac{(k^1_i)^{\delta}-(k^2_i)^{\delta}}{(k^1_i)^{\delta}}
  \frac{(k^2_{i'})^{\delta}-(k^1_{i'})^{\delta}}{\lambda(T)}
  \times\frac{(k^1_i)^{\delta}}{l_{\delta}}\\
  && \quad = \quad \frac{1}{l_\delta}\sum_{i'\in R_2} 
  ((k^2_{i'})^{\delta}-(k^1_{i'})^{\delta}).
\end{eqnarray}
Since, for any concave function $f: \R_+ \rightarrow \R$
and any $s\in\N\setminus\{0\}$,
$(z_1,\dots,z_s)\in\R_+^s\mapsto f(\sum_i z_i) - \sum_i f(z_i)$
is non-increasing in all its $s$ variables (as a consequence
of the slope inequalities), by concavity of $z\mapsto z^\delta$
and using the fact that, for all $i'\in R_2$,
$k^2_{i'} > k^1_{i'}$ we get
\begin{eqnarray}
  P(i^1 \in R_1, i^2 \in R_2)
  & \geq & 
  \frac{1}{l_\delta}\left[
    \bigg(\sum_{i'\in R_2} k^2_{i'}\bigg)^{\delta}
    - \bigg(\sum_{i'\in R_2} k^1_{i'}\bigg)^{\delta}
  \right]\\
  & = &
  \frac{1}{l_\delta}\left[
    \bigg(\rho_t + \sum_{i'\in R_2} k^1_{i'}\bigg)^{\delta}
    - \bigg(\sum_{i'\in R_2} k^1_{i'}\bigg)^{\delta}
  \right].
\end{eqnarray}
Using the same property of concave functions on $\R^+$
(with $s = 2$) and the fact that 
$\sum_{i'\in R_2} k^1_{i'} \leq k - \rho_t$,
then using once again the concavity
of $z\mapsto z^\delta$,
we write
\begin{equation}
  P(i^1 \in R_1, i^2 \in R_2)
  \geq  
  \frac{1}{l_\delta}\left[
    k^{\delta}
    - (k-\rho_t)^{\delta}
  \right]
  \geq \frac{\delta k^{\delta-1}\rho_t}{l_\delta}
\end{equation}
so that, by the previous proposition,
\begin{equation}
  P(\rho_{t+1} = \rho_t - 1 | X^1_t, X^2_t)
  \geq \frac{\delta k^{\delta-1} \rho_t}{ml_{\delta}}.
\end{equation}
\qed

\medskip\par\noindent
As a consequence
\begin{equation}
  E(\rho_{t+1} | X^1_t, X^2_t) 
  \leq \left(1 - \frac{\delta}{km(k \wedge m)^{1-\delta}}\right)\rho_t
\end{equation}
and, by Proposition~\ref{prpdist}, 
this finally proves Theorem~\ref{thm2}.

\begin{appendix}
\section{Simulation appendix}
In this appendix we report on simulations
for the temporal evolution of $d(t)$
defined in~(\ref{defd})
as well as the dependency 
in $\beta$ of the mixing time $t_\epsilon$
defined in~(\ref{deft})
with $\epsilon=0.1$ 
in the context of 
the Fermi statistics $\nu$
associated with $k=50$ particles
among $m=20$ energy levels
with two different profiles.
\begin{itemize}
\item
  In the first case we chose $v_j=j/20$
  and $n_j = 2^j$ for $1\leq j\leq m$
  so that $n = \sum_j 2^j = 2^{21} - 2$. See Figure~\ref{fig1}. 
  \begin{figure}[H]
   \includegraphics[width=0.48\textwidth]{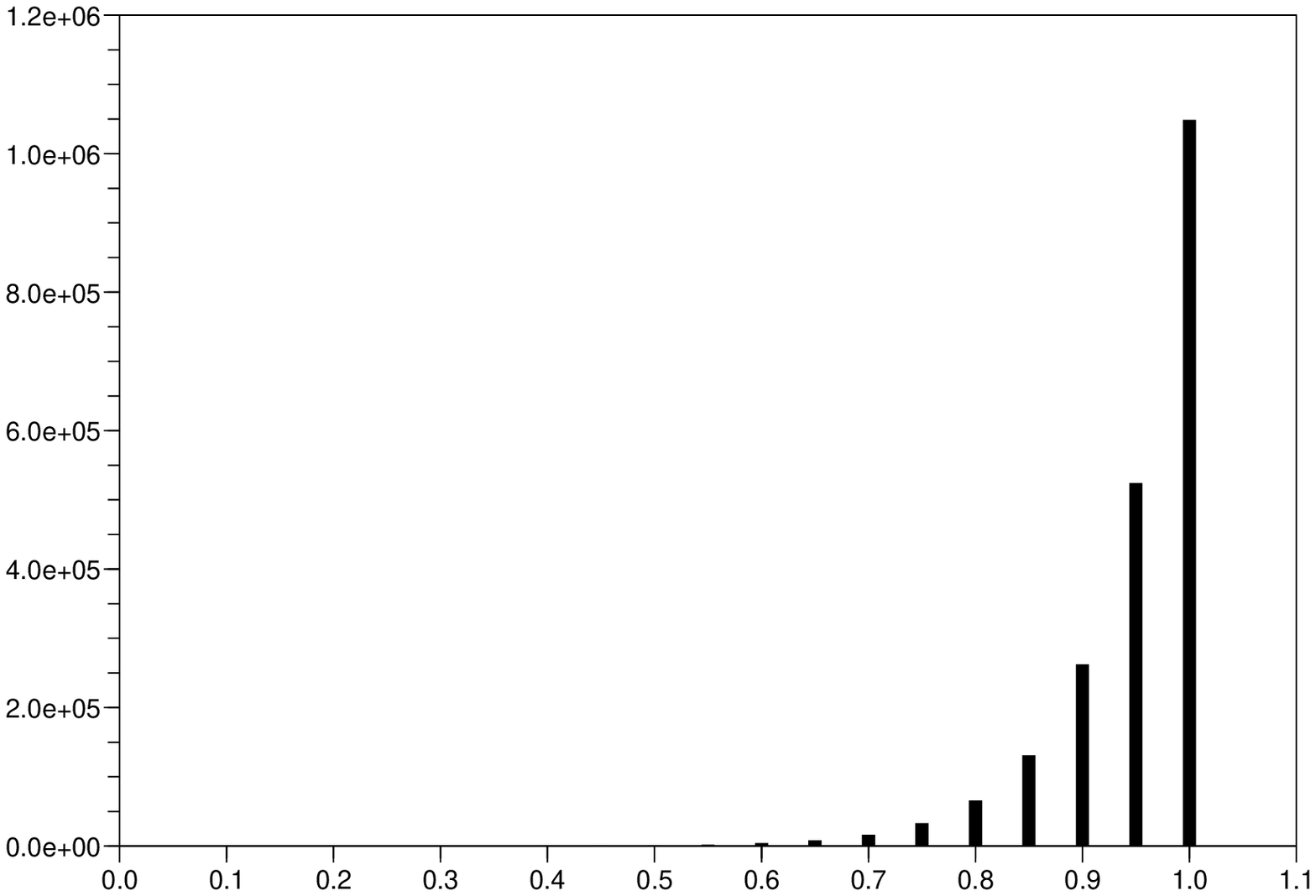}
   \includegraphics[width=0.48\textwidth]{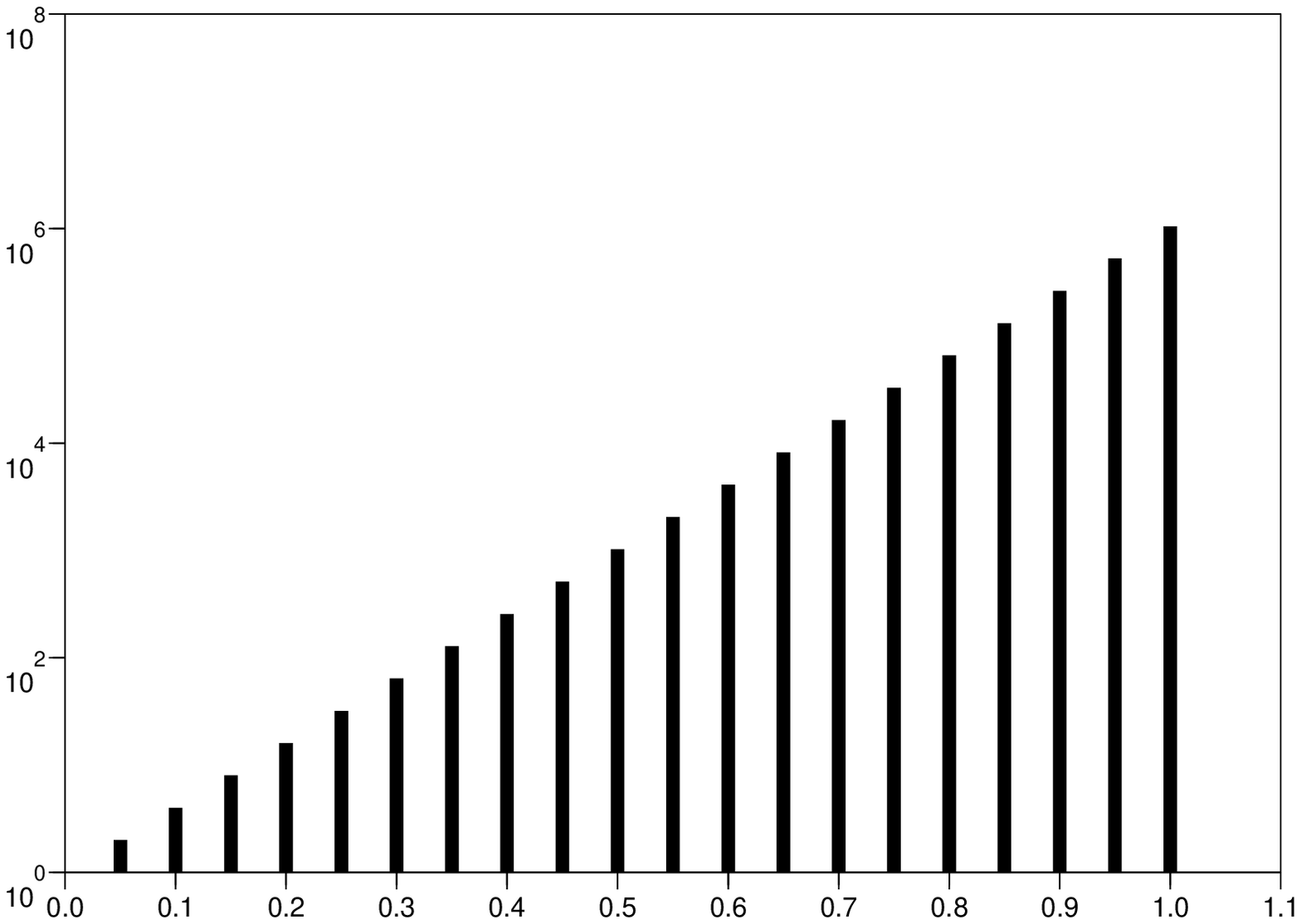}
   \caption{Degeneracies of the energy levels for case 1, in semi-logarithmic
     representation for the second picture.}
   \label{fig1}
 \end{figure}
\item
  In the second case we chose the same $v_j$'s, 
  the same volume $n = 2^{21} - 2$
  and we chose the degeneracies according
  two a multinomial law with parameters
  $n$, $1/20$,~\dots,~$1/20$. See Figure~\ref{fig2}.
  \begin{figure}[H]
   \includegraphics[width=0.48\textwidth]{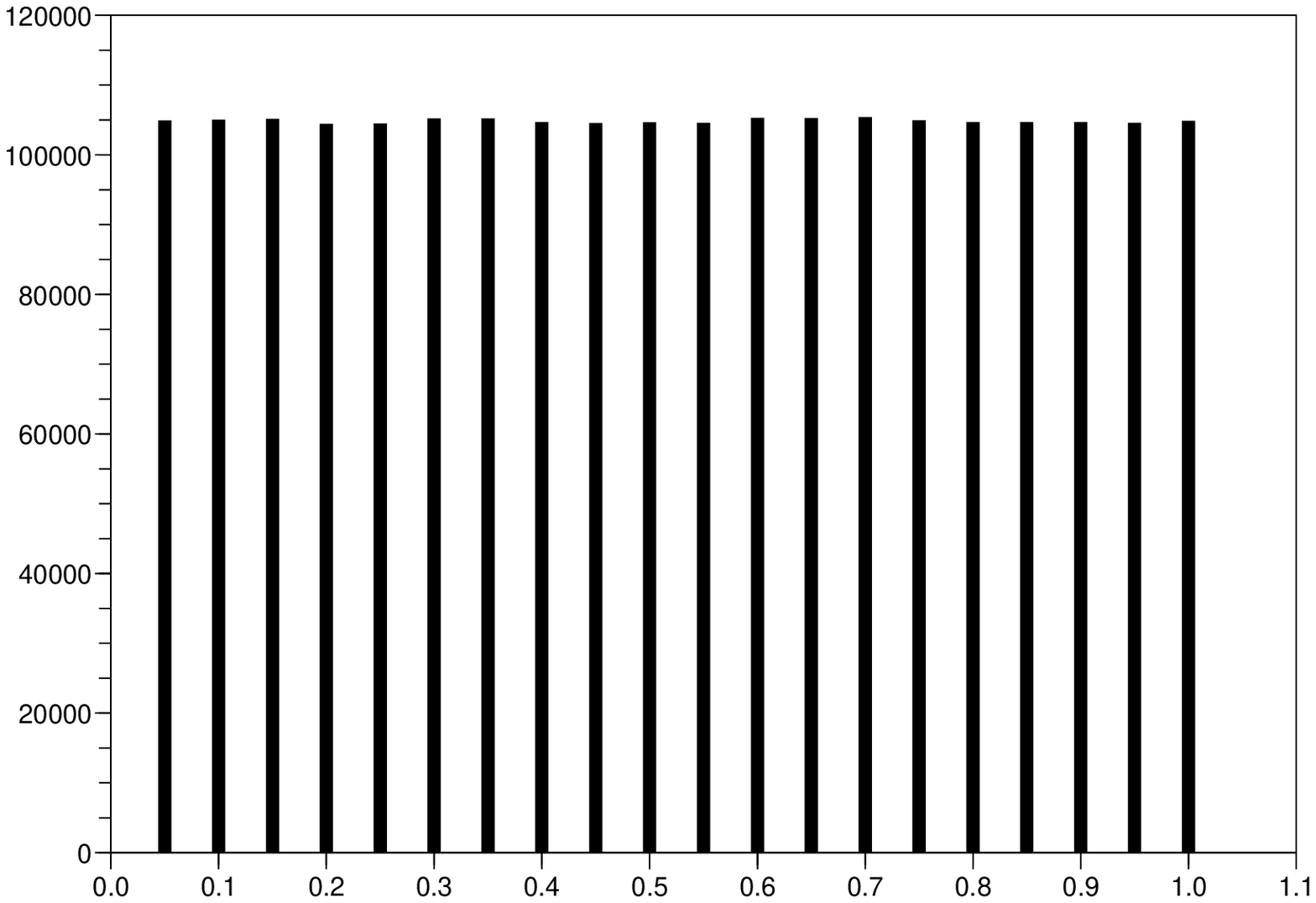}
   \includegraphics[width=0.48\textwidth]{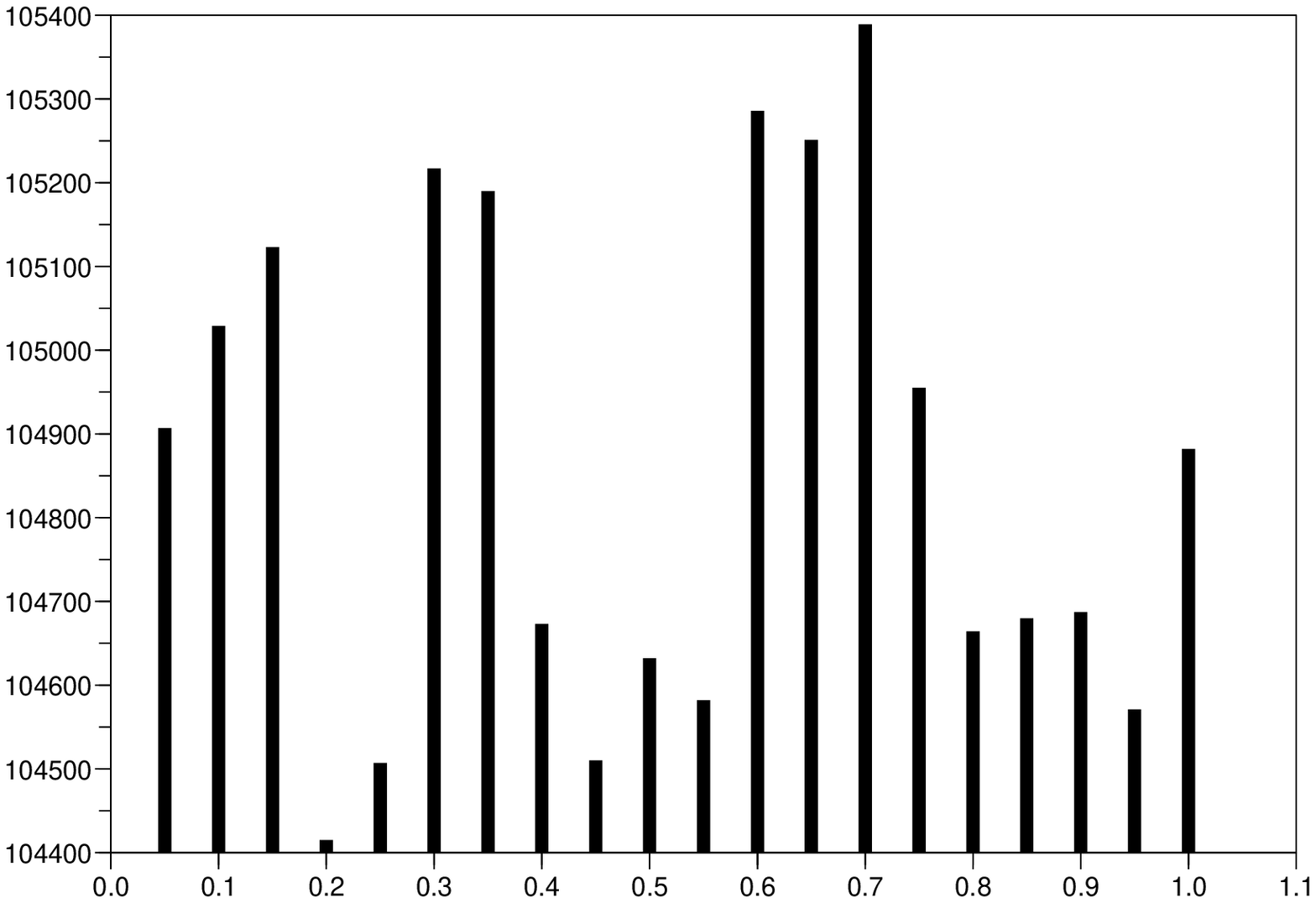}
   \caption{Degeneracies of the energy levels for case 2,
     with shifted and rescaled axis for the second picture.}
   \label{fig2}
 \end{figure}
\end{itemize}

We were then interested in
\begin{equation}
d(t)
= \max_{\eta\in{\cal X}_{k,m}}\|p^t(\eta_0,\cdot) - \nu\|_{TV} 
= \max_{\eta\in{\cal X}_{k,m}}\frac{1}{2}\sum_{\eta\in{\cal X}_{k,m}}|p^t(\eta_0,\eta) - \nu(\eta)| 
\label{3pbs}
\end{equation}
for $t\geq 0$.
The cardinality of ${\cal X}_{k,m}$
being very large (equal to $\binom{69}{19}\simeq 5\times 10^{16}$ in case 2),
one faces three difficulties with such a formula:
\begin{itemize}
\item[i)]
  One cannot compute $\nu(\eta)$ and $p^t(\eta_0,\eta)$ for given
  $\eta$ and $\eta_0$.
\item[ii)]
  The sum in~(\ref{3pbs}) contained too many terms to be computed.
\item[iii)]
  It is not possible to try all the possible $\eta_0$ 
  before taking the maximum in~(\ref{3pbs}).
\end{itemize}
We addressed the first difficulty
by replacing $p^t(\eta_0,\cdot)$ by
the empirical measure 
\begin{equation}
\mu_t^N:=\frac{1}{N}\sum_{l=1}^{N}\delta_{X_t^l}
\end{equation}
computed with a large number $N$ of simulated
independent copies $(X_t^l)_{t\in{\mathbb N}}$ of our Markov chain
starting from $\eta_0$ and, using Theorem~\ref{thm1},
by replacing $\nu$ by 
\begin{equation}
\nu^N:=\frac{1}{2T}\sum_{t=T+1}^{2T}{\mu'_t}^N+{\mu''_t}^N
\end{equation}
with $T:=\lfloor km\ln(k/\epsilon)\rfloor$
and where ${\mu'_t}^N$
and ${\mu''_t}^N$
are independently computed like 
$\mu_t^N$,
starting from two different configurations.
We addressed the second difficulty 
by looking at a coarse grained version
$\bar{\cal X}_{k,m}$ of ${\cal X}_{k,m}$
through the free entropy $\phi$.
More precisely, with $I^N = [\phi_-,\phi_+]$
the smallest interval containing the support
of $\phi(\nu^N)= \nu^N\circ\phi^{-1}$
we divided $I^N$ into $M$ intervals
of length $\Delta = (\phi_+-\phi_-)/M$,
we extended this partition  of $I^N$
into a partition of ${\mathbb R}$
in intervals of length $\Delta$
and grouped in a same class the configurations
$\eta\in{\cal X}_{k,m}$ that fall into a same interval.
(For possible cases where $\phi_- = \phi_+$,
we used a partition of ${\mathbb R}$ into
intervals of length $10^{-6}$.)
With $\bar\phi$ the canonical projection
from ${\cal X}_{k,m}$ to $\bar{\cal X}_{k,m}$
we could then compute
$\|\bar\phi(\mu_t^N) - \bar\phi(\nu^N)\|_{TV}$
rather than $\|\mu_t^N -\nu^N\|_{TV}$.
Note that 
\begin{equation}
\|\bar\phi(\mu_t^N) - \bar\phi(\nu^N)\|_{TV}
\leq\|\mu_t^N -\nu^N\|_{TV},
\end{equation}
and we also have, for large $M$ and as soon
as $\phi$ is an injection from ${\cal X}_{k,m}$ to ${\mathbb R}$,
\begin{equation}
\|\mu_t^N - \nu^N\|_{TV} 
= \|\bar\phi(\mu_t^N)-\bar\phi(\nu^N)\|_{TV} + o(1).
\end{equation}
As far as the third difficulty was concerned
we underestimated
\begin{equation}
\bar d^N(t)
= \max_{\eta_0\in{\cal X}_{k,m}}
\|\bar\phi(\mu_t^N) - \bar\phi(\nu^N)\|_{TV}
\end{equation}
by making a guess on a particular configuration
$\eta_0$ for which $\|\bar\phi(\mu_t^N) - \bar\phi(\nu^N)\|_{TV}$
could be of the order of $\bar d^N(t)$.
We simply took an $\eta_0$ for which we could have say,
for all $\beta$,  that it was ``very far from typical
configurations''.
We first took $\eta_0 = (0,\dots,0,k)$ since typical configurations
are concentrated on the low energy levels when $\beta$ is large
and have their occupation numbers $k_j$ distributed like the
degeneracies $n_j$ when $\beta$ is small.
Since this last point tends to show that in case~1
our specific $\eta_0$ is ``not so far'' from typical equilibrium
configurations (we will come back on this point later),
we also considered a different $\eta_0$,
that one for which the high energy levels are empty
and the low energy levels are saturated 
or as close to saturation as possible: 
\begin{equation}
  \eta_0=\left\{
    \begin{array}{ll}
      (2, 4, 8, 16, 20, 0, \dots, 0) & \mbox{in case 1,}\\
      (50, 0, \dots, 0) & \mbox{in case 2.}
    \end{array}
  \right.
\end{equation}
Then we simply took the maximum 
of the two quantities $\|\bar\phi(\mu_t^N) - \bar\phi(\nu^N)\|_{TV}$
computed, in each case, with these two specific choices.

\begin{figure}[H]
  \includegraphics[width=0.48\textwidth]{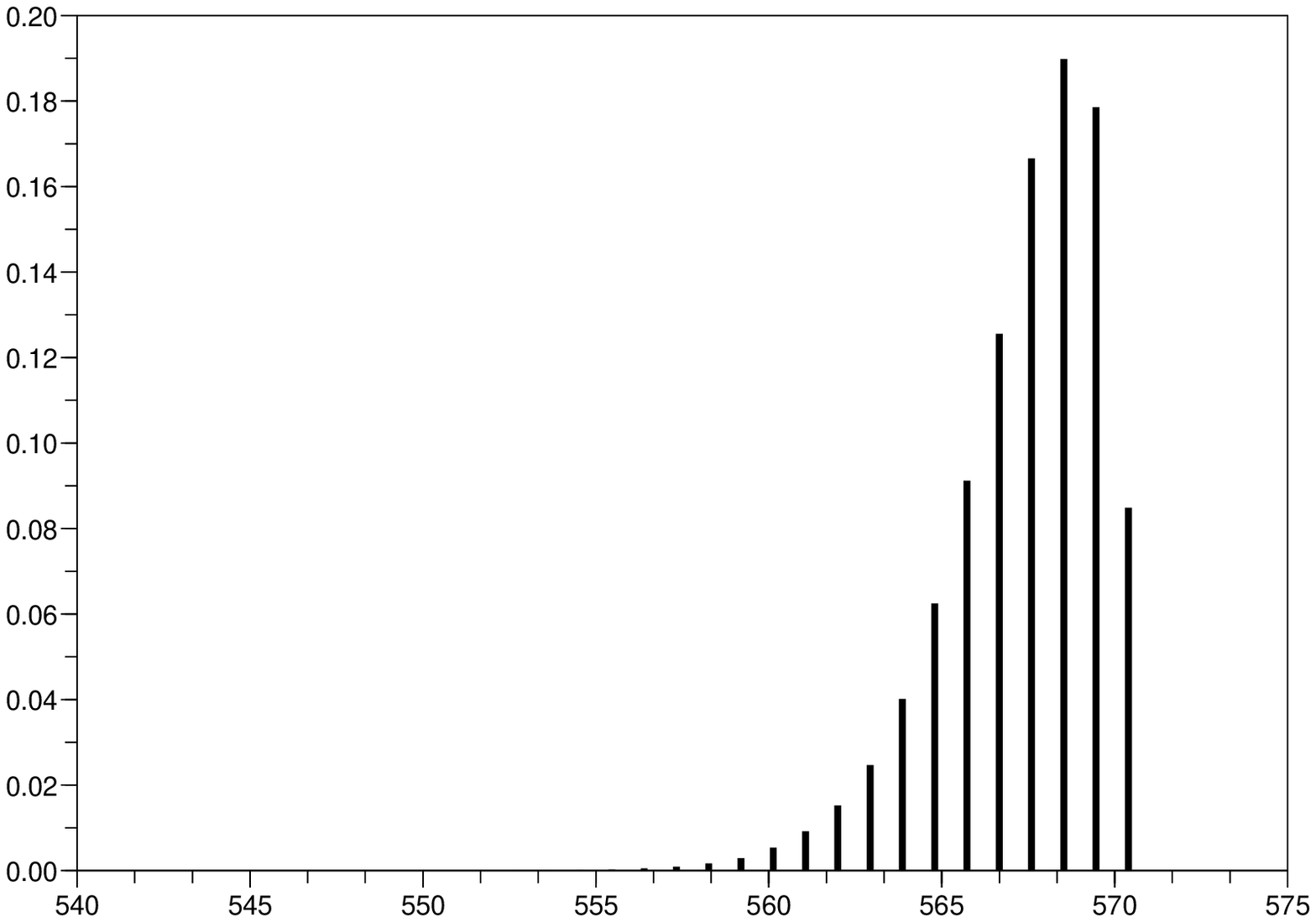}
  \includegraphics[width=0.48\textwidth]{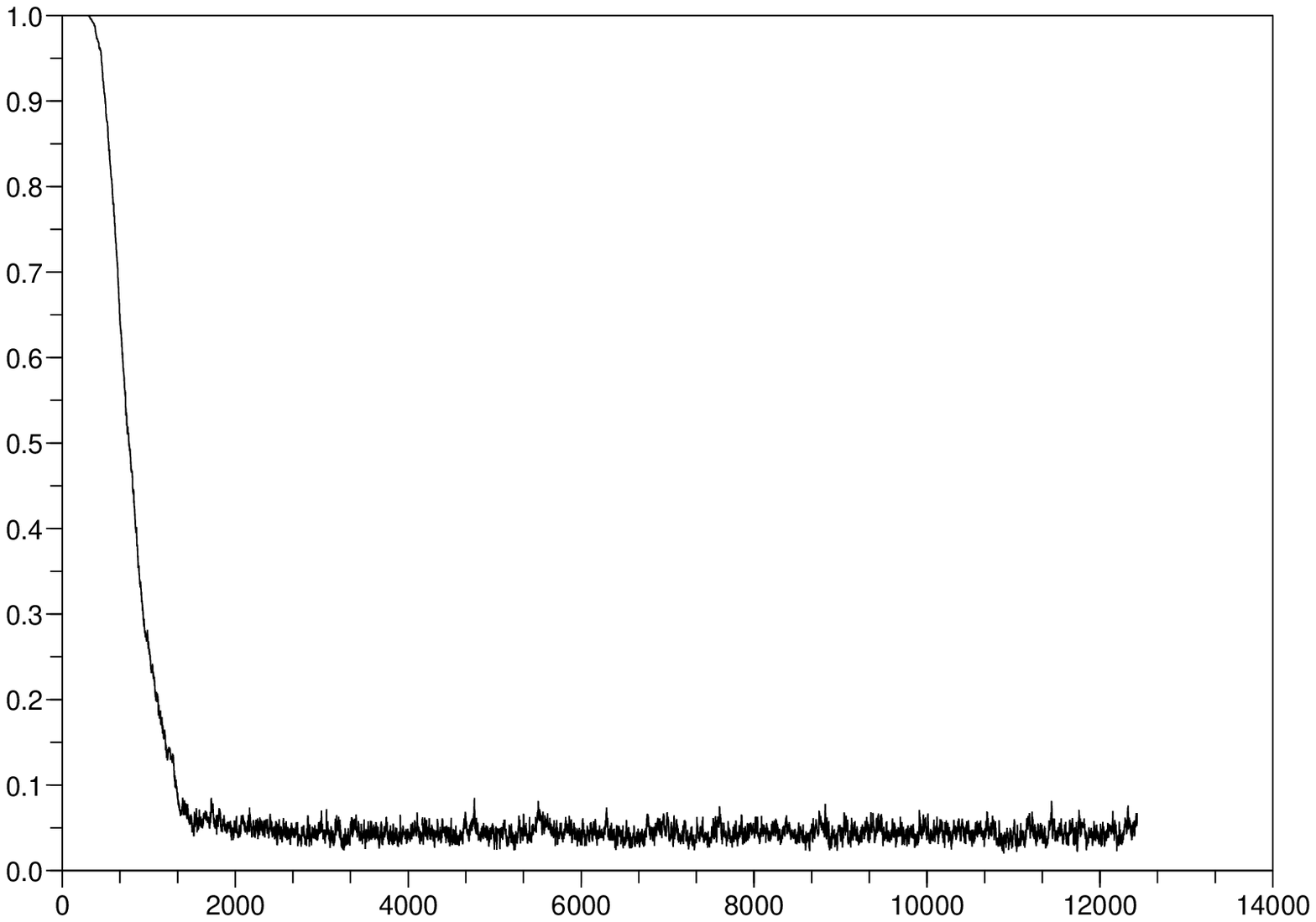}
  \caption{Probability distribution function of $\bar\phi(\nu^N)$
    and temporal evolution of $\bar d_0^N(t)$
    for case 1 and $\beta = 0$.}
  \label{fig3.1}
\end{figure}

\begin{figure}[H]
  \includegraphics[width=0.48\textwidth]{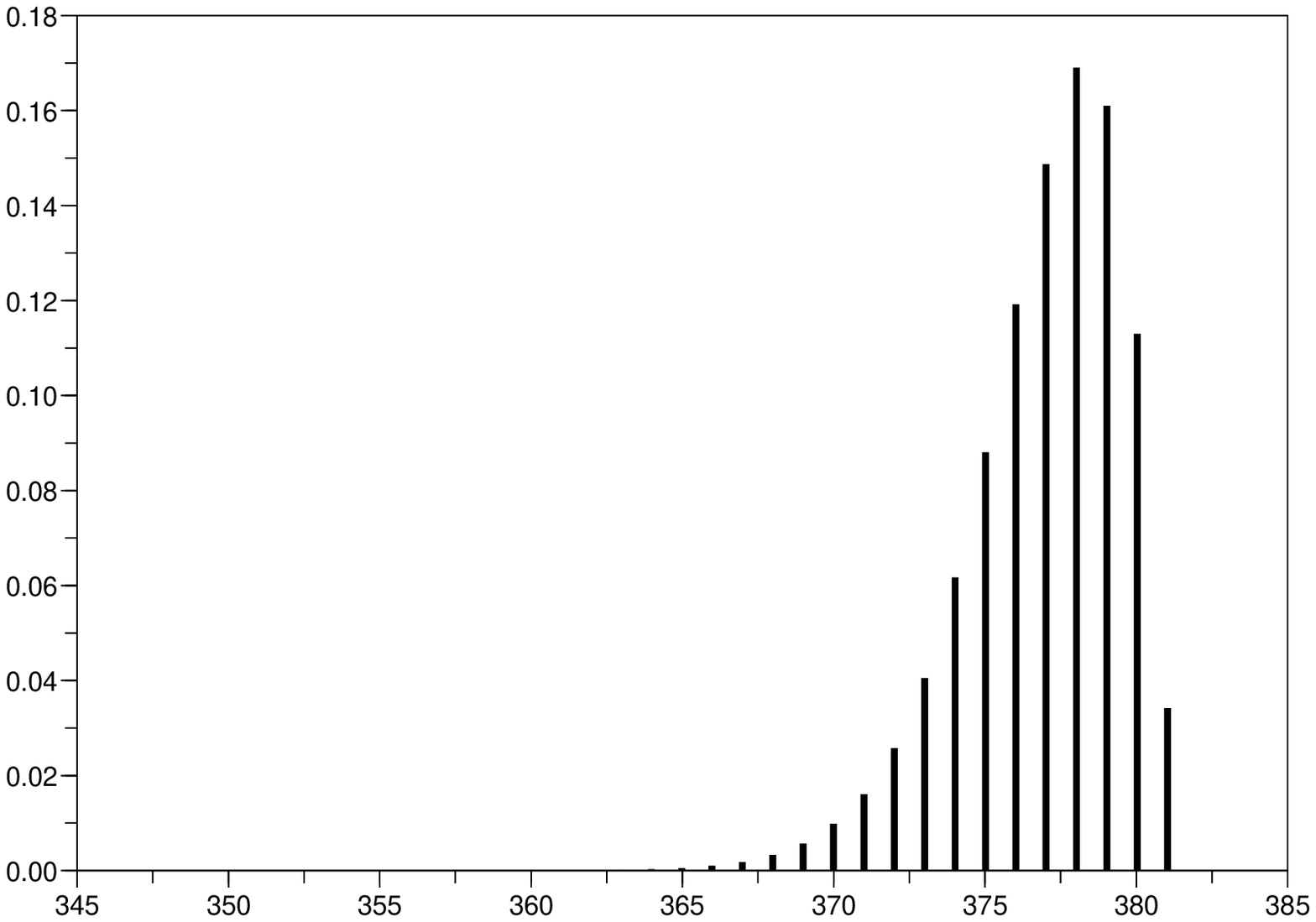}
  \includegraphics[width=0.48\textwidth]{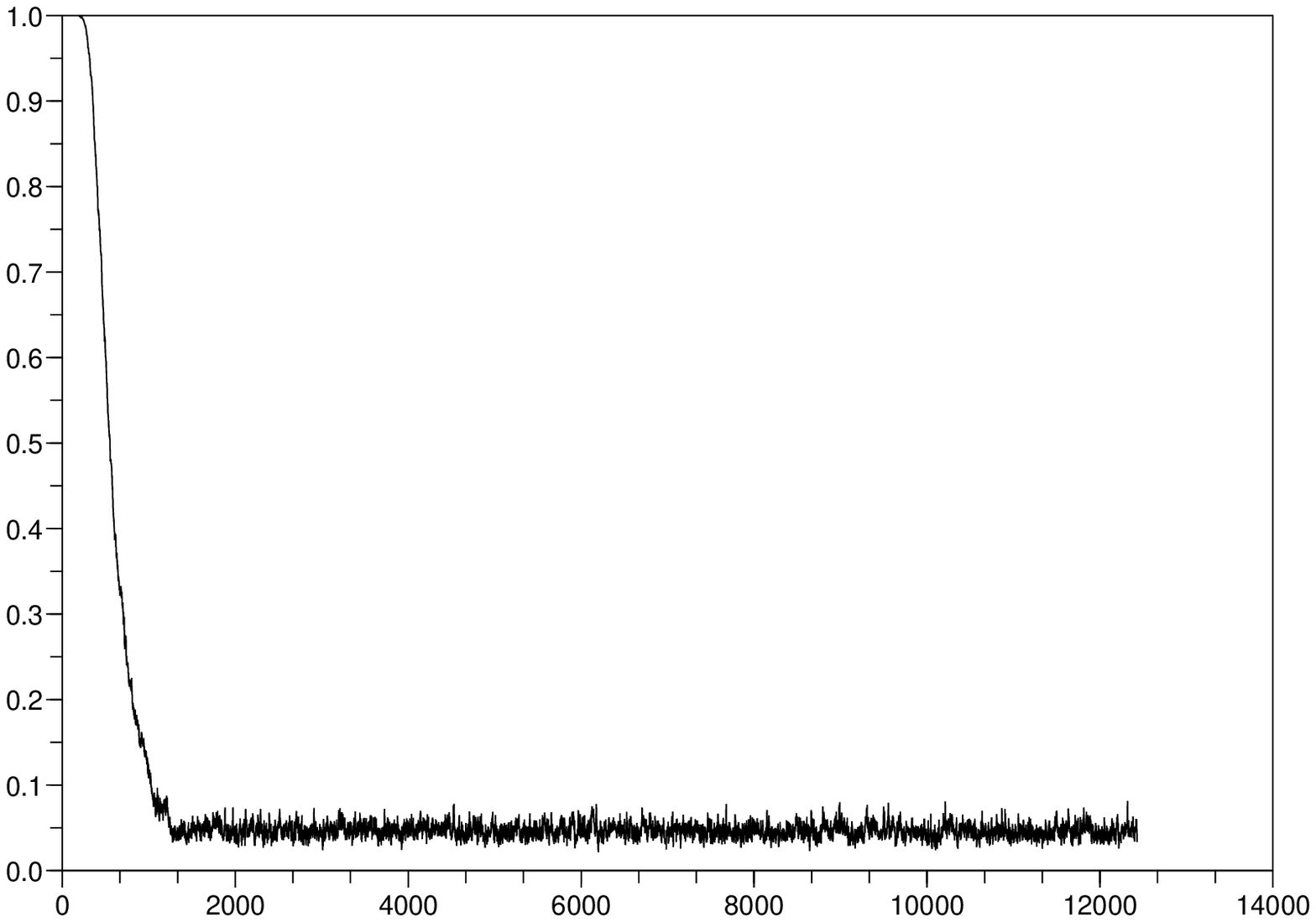}
  \caption{Probability distribution function of $\bar\phi(\nu^N)$
    and temporal evolution of $\bar d_0^N(t)$
    for case 1 and $\beta = 4$.}
  \label{fig3.2}
\end{figure}

\begin{figure}[H]
  \includegraphics[width=0.48\textwidth]{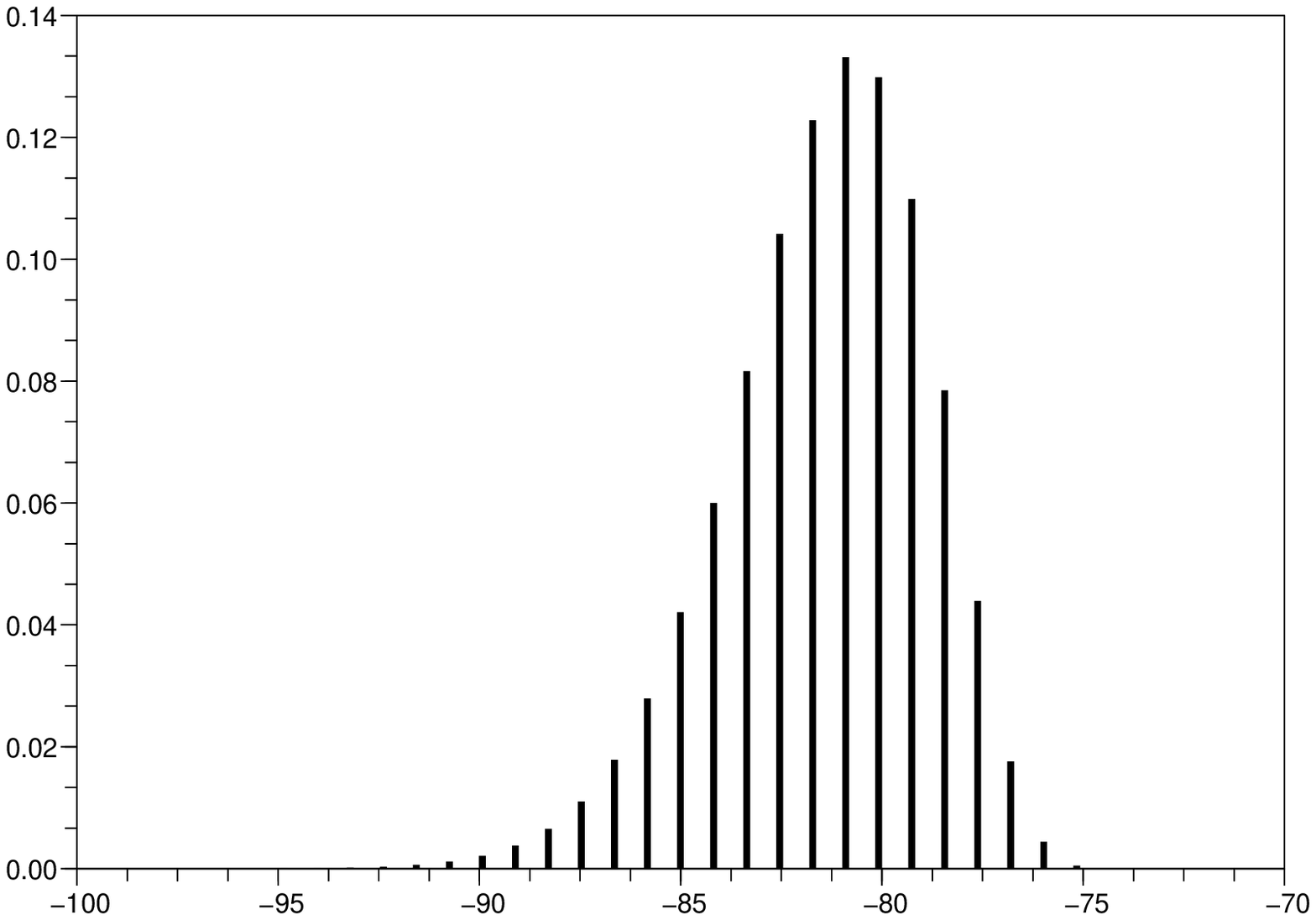}
  \includegraphics[width=0.48\textwidth]{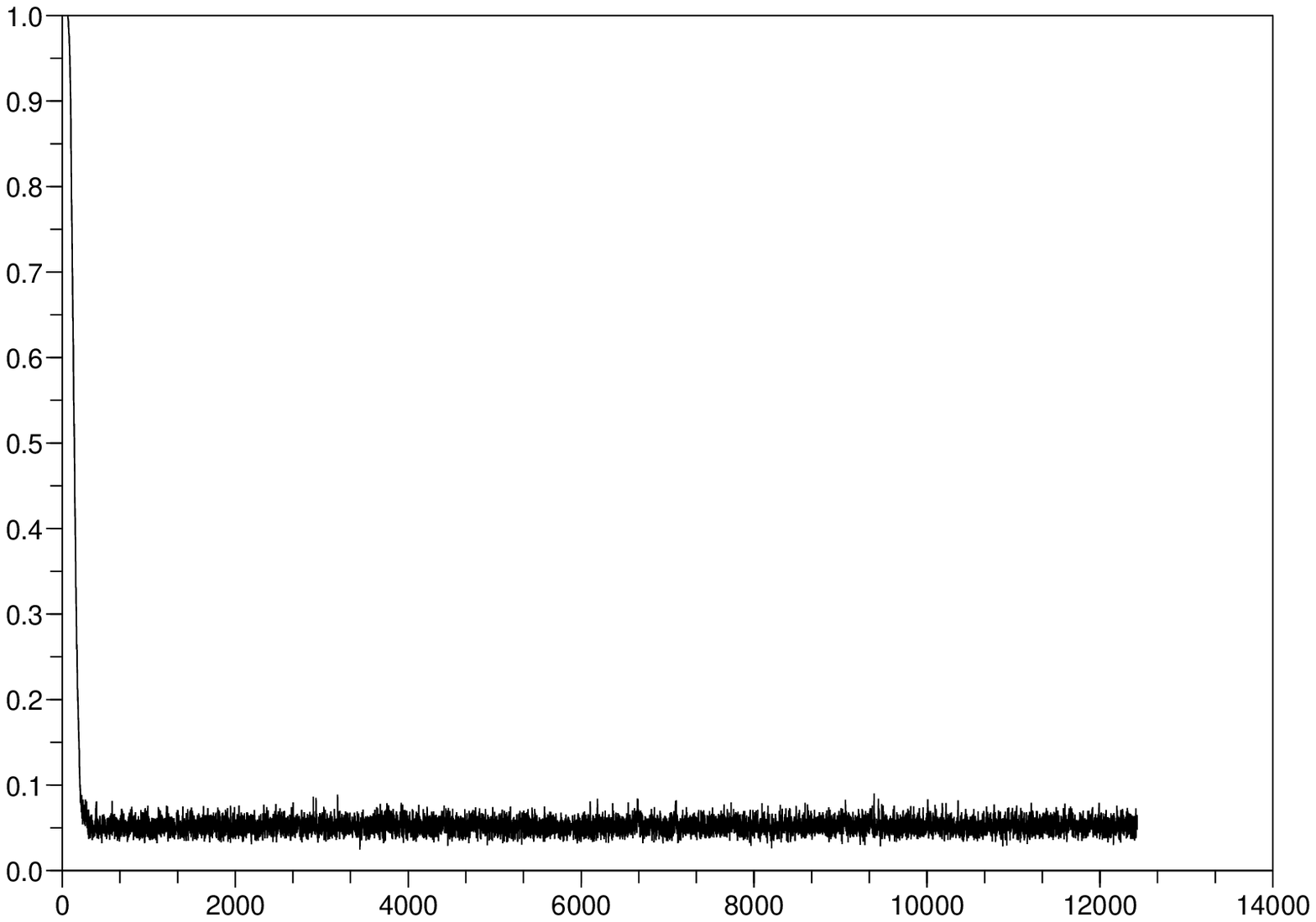}
  \caption{Probability distribution function of $\bar\phi(\nu^N)$
    and temporal evolution of $\bar d_0^N(t)$ 
    for case 1 and $\beta = 16$.}
  \label{fig3.3}
\end{figure}

\begin{figure}[H]
  \includegraphics[width=0.48\textwidth]{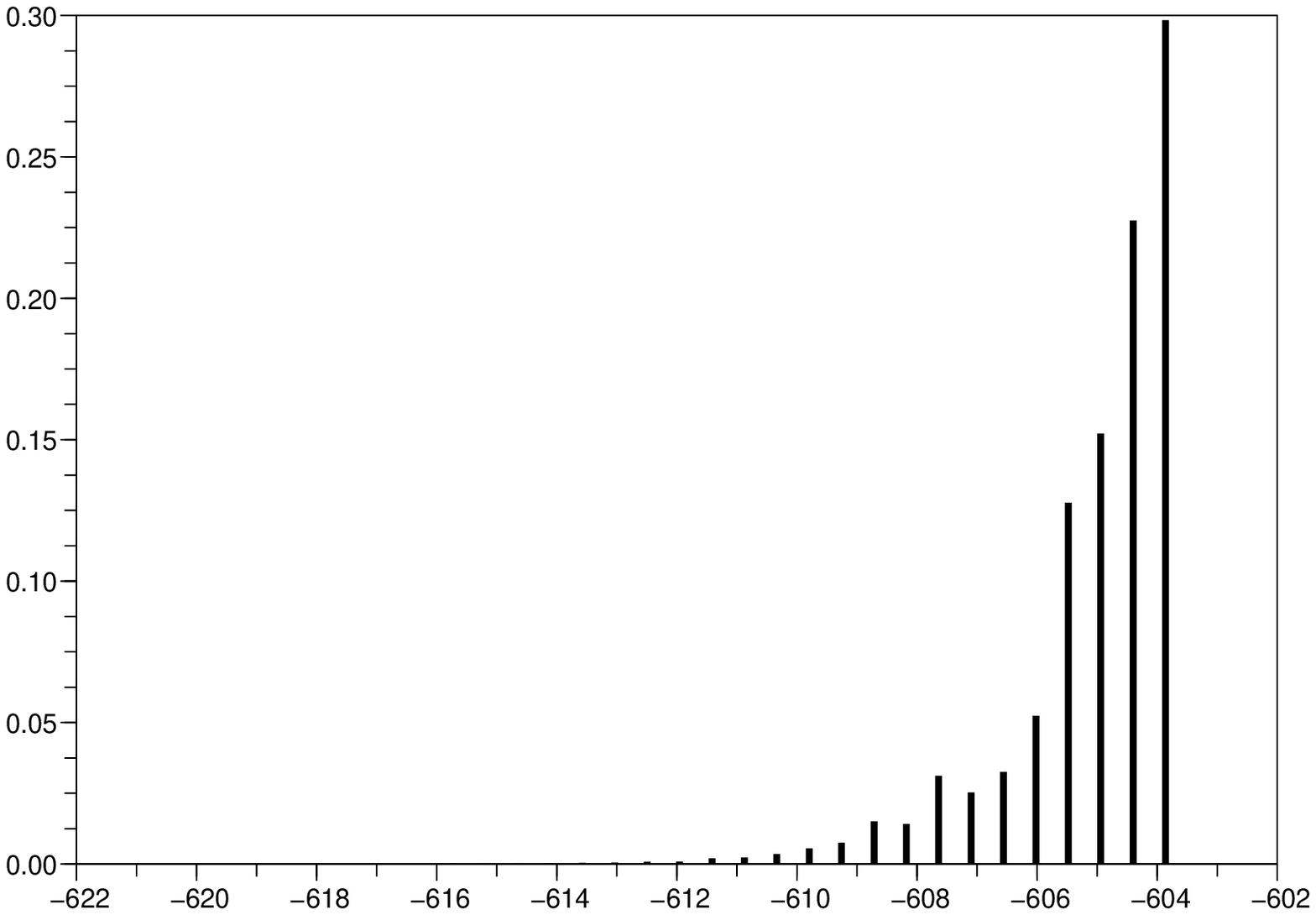}
  \includegraphics[width=0.48\textwidth]{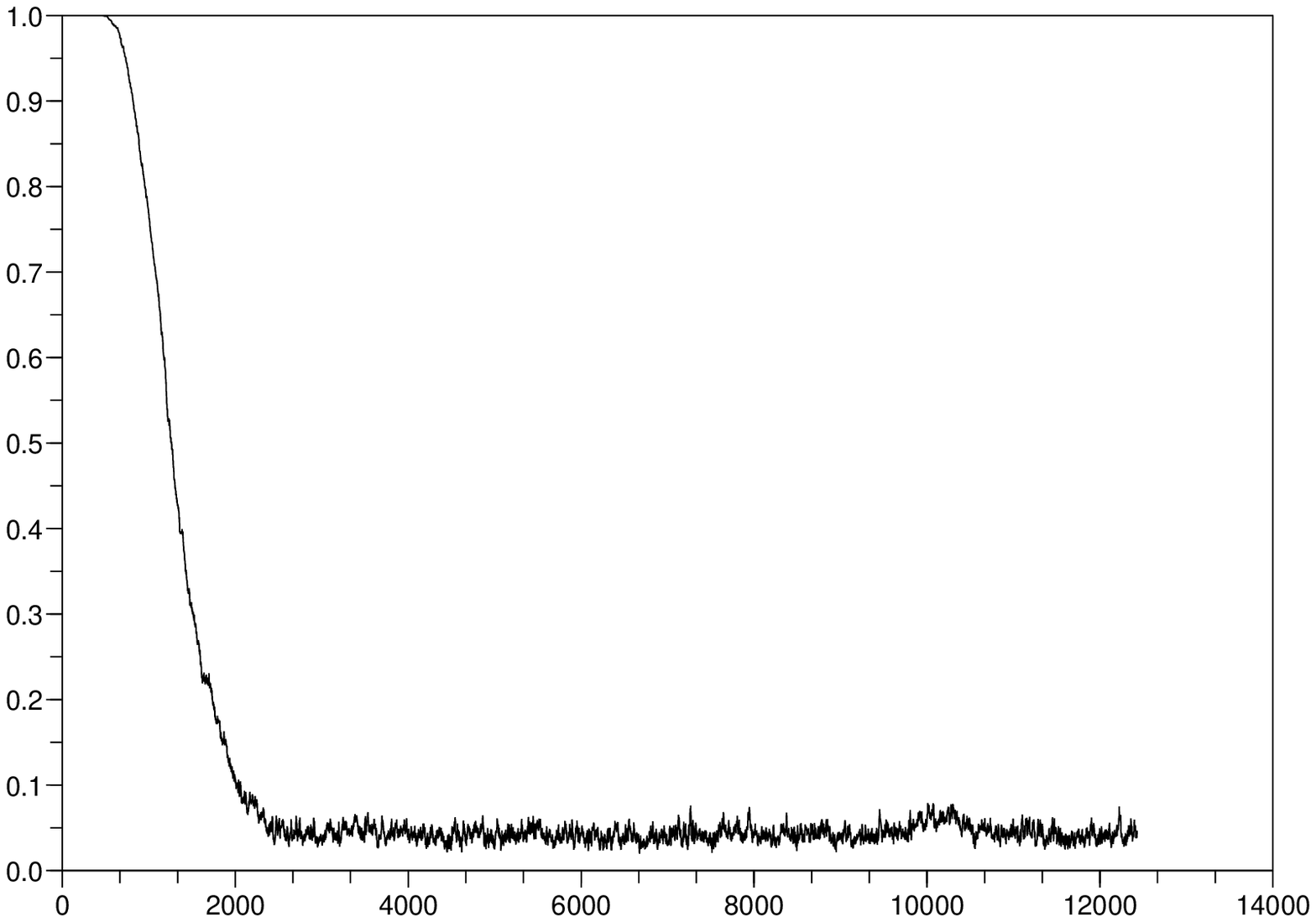}
  \caption{Probability distribution function of $\bar\phi(\nu^N)$
    and temporal evolution of $\bar d_0^N(t)$
    for case 1 and $\beta = 63$.}
  \label{fig3.4}
\end{figure}

\begin{figure}[H]
  \includegraphics[width=0.48\textwidth]{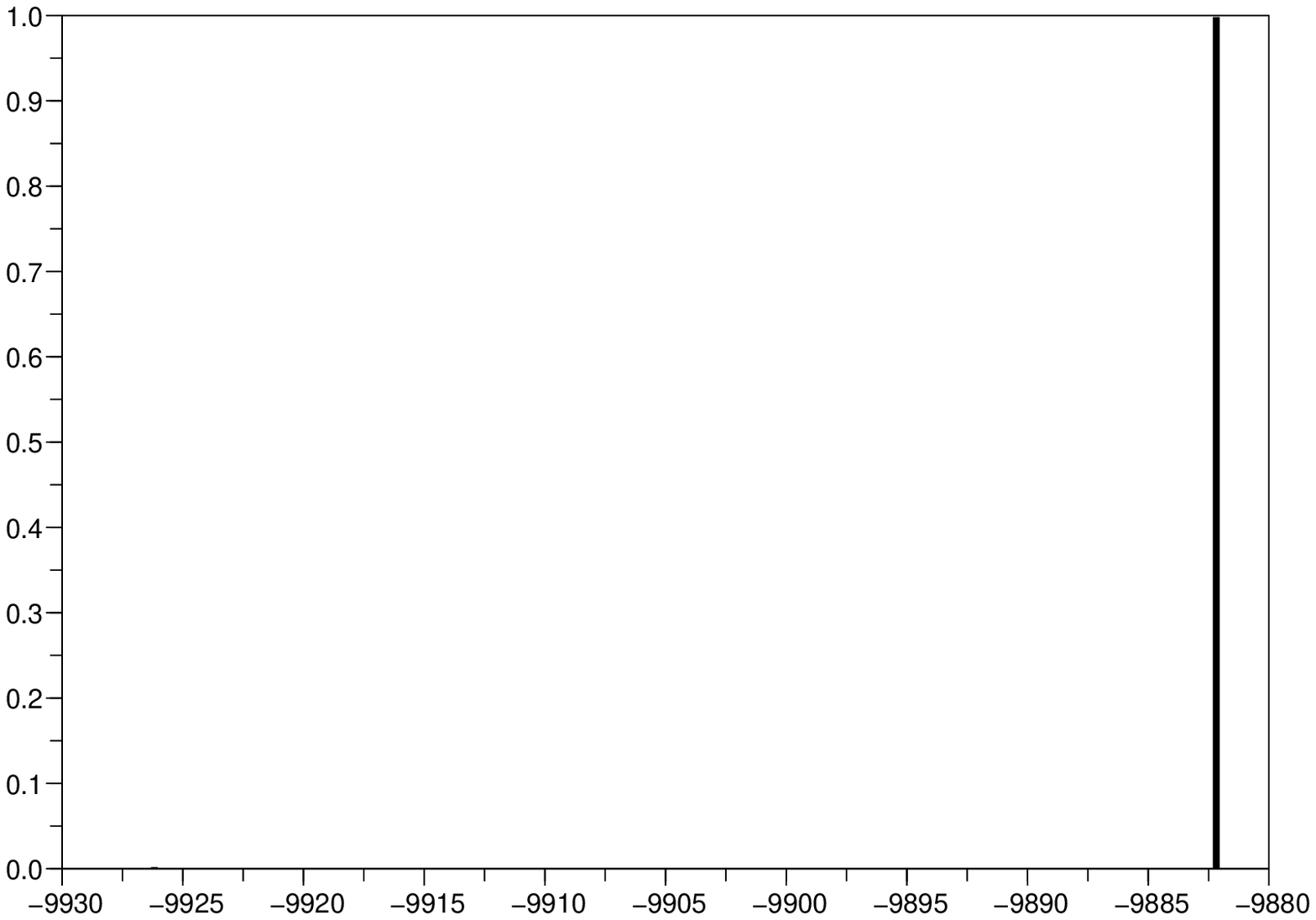}
  \includegraphics[width=0.48\textwidth]{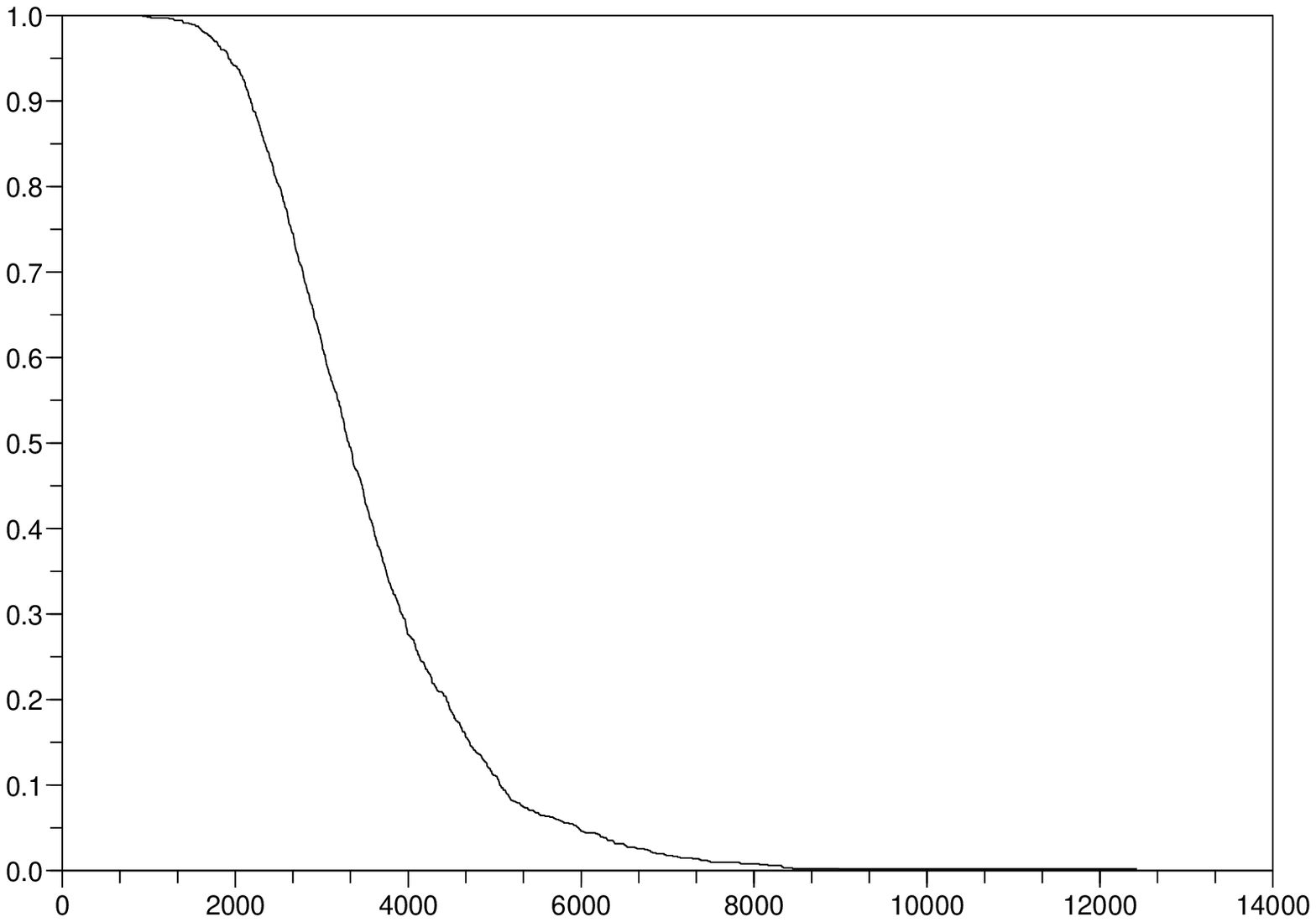}
  \caption{Probability distribution function of $\bar\phi(\nu^N)$
    and temporal evolution of $\bar d_0^N(t)$
    for case 1 and $\beta = 1000$.}
  \label{fig3.5}
\end{figure}

\begin{figure}[H]
  \includegraphics[width=0.48\textwidth]{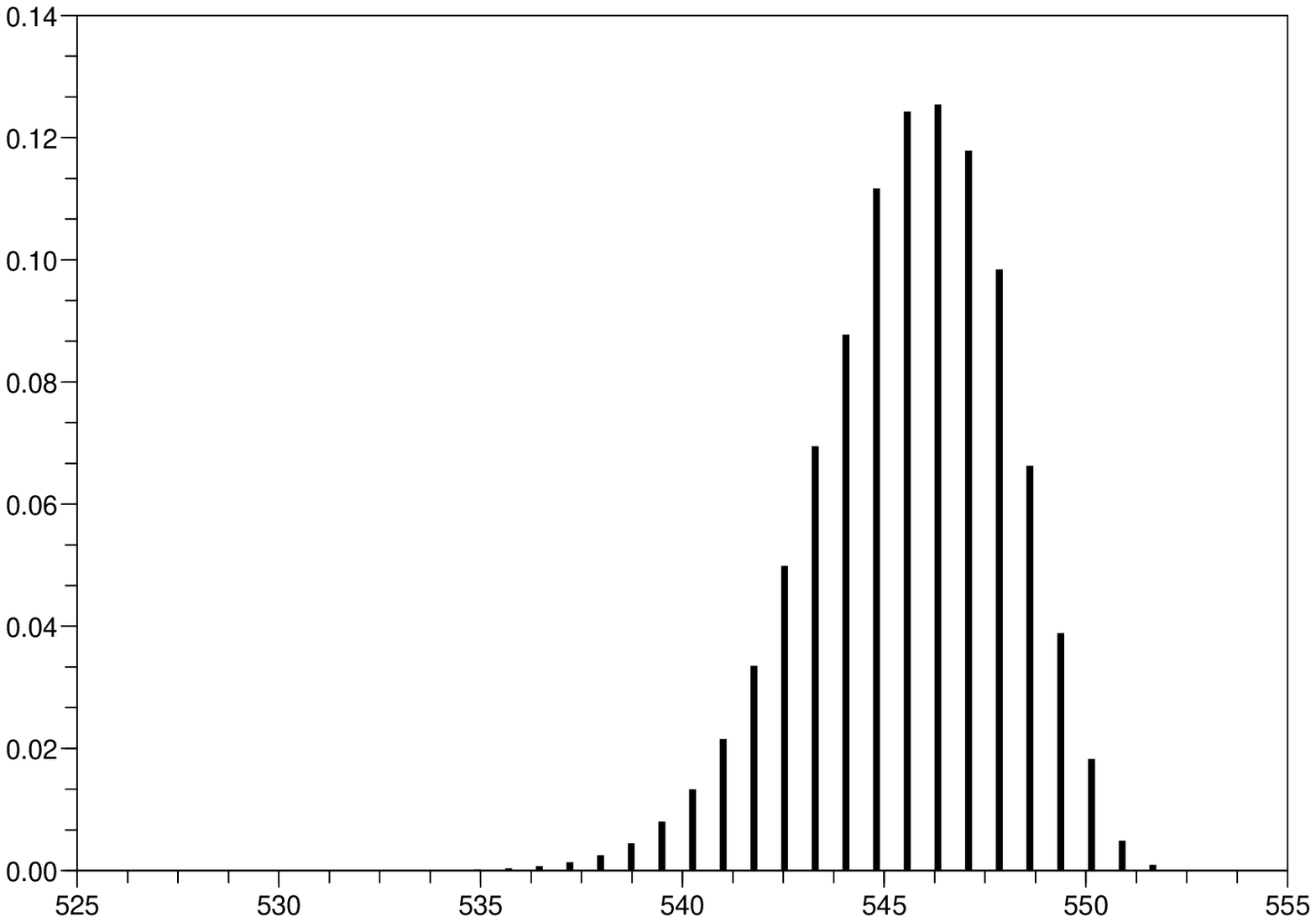}
  \includegraphics[width=0.48\textwidth]{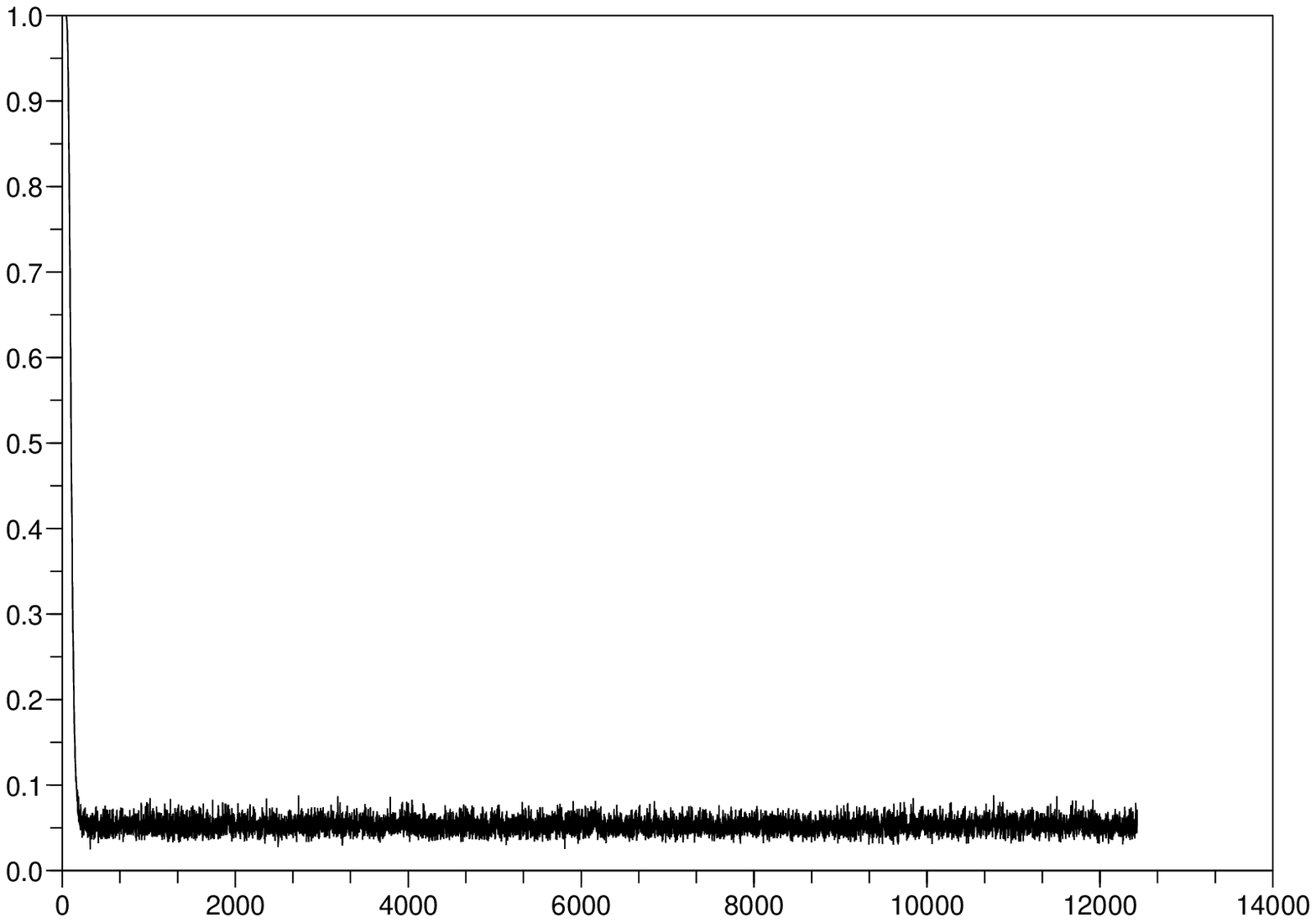}
  \caption{Probability distribution function of $\bar\phi(\nu^N)$
    and temporal evolution of $\bar d_0^N(t)$
    for case 2 and $\beta = 0$.}
  \label{fig4.1}
\end{figure}

\begin{figure}[H]
  \includegraphics[width=0.48\textwidth]{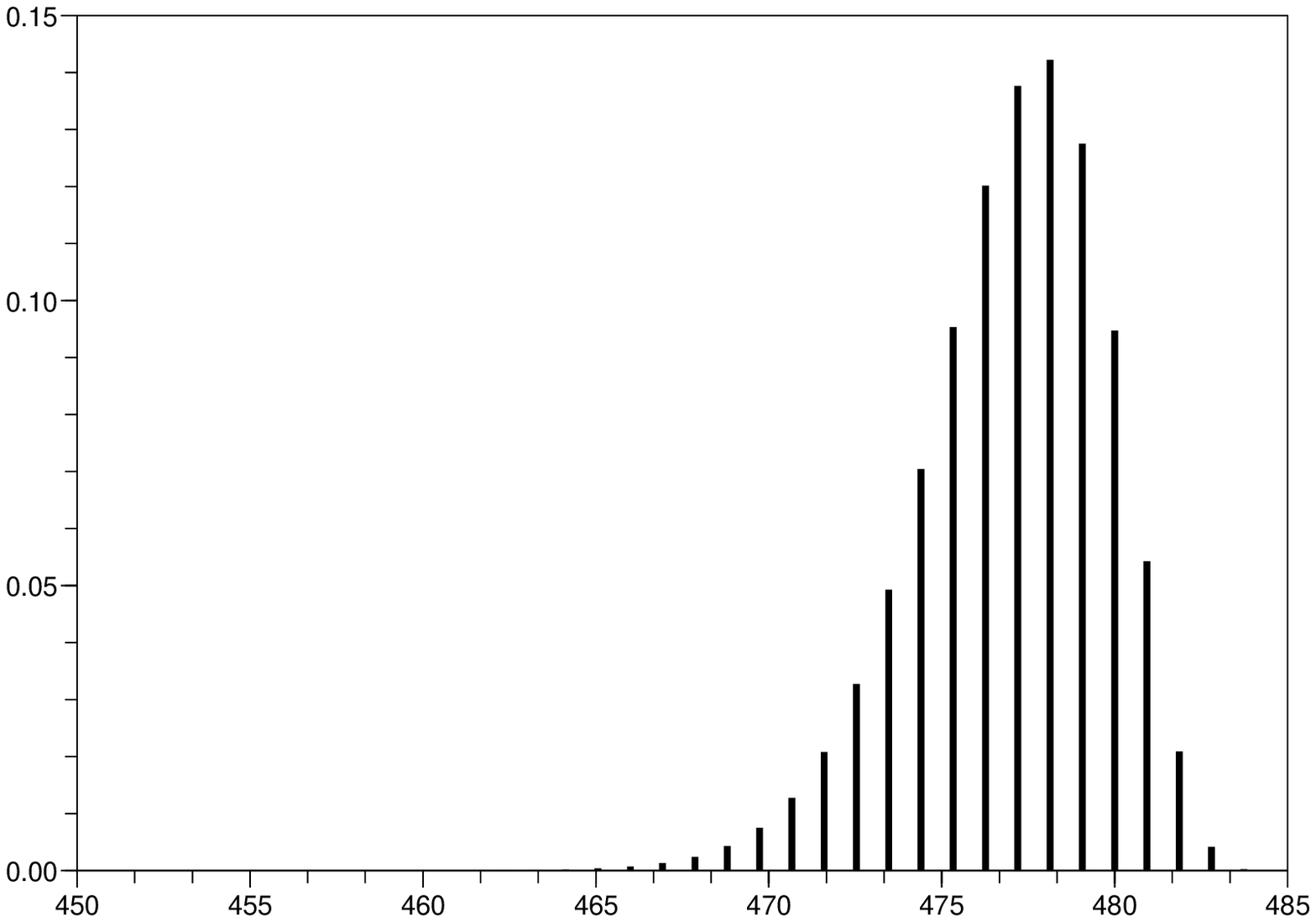}
  \includegraphics[width=0.48\textwidth]{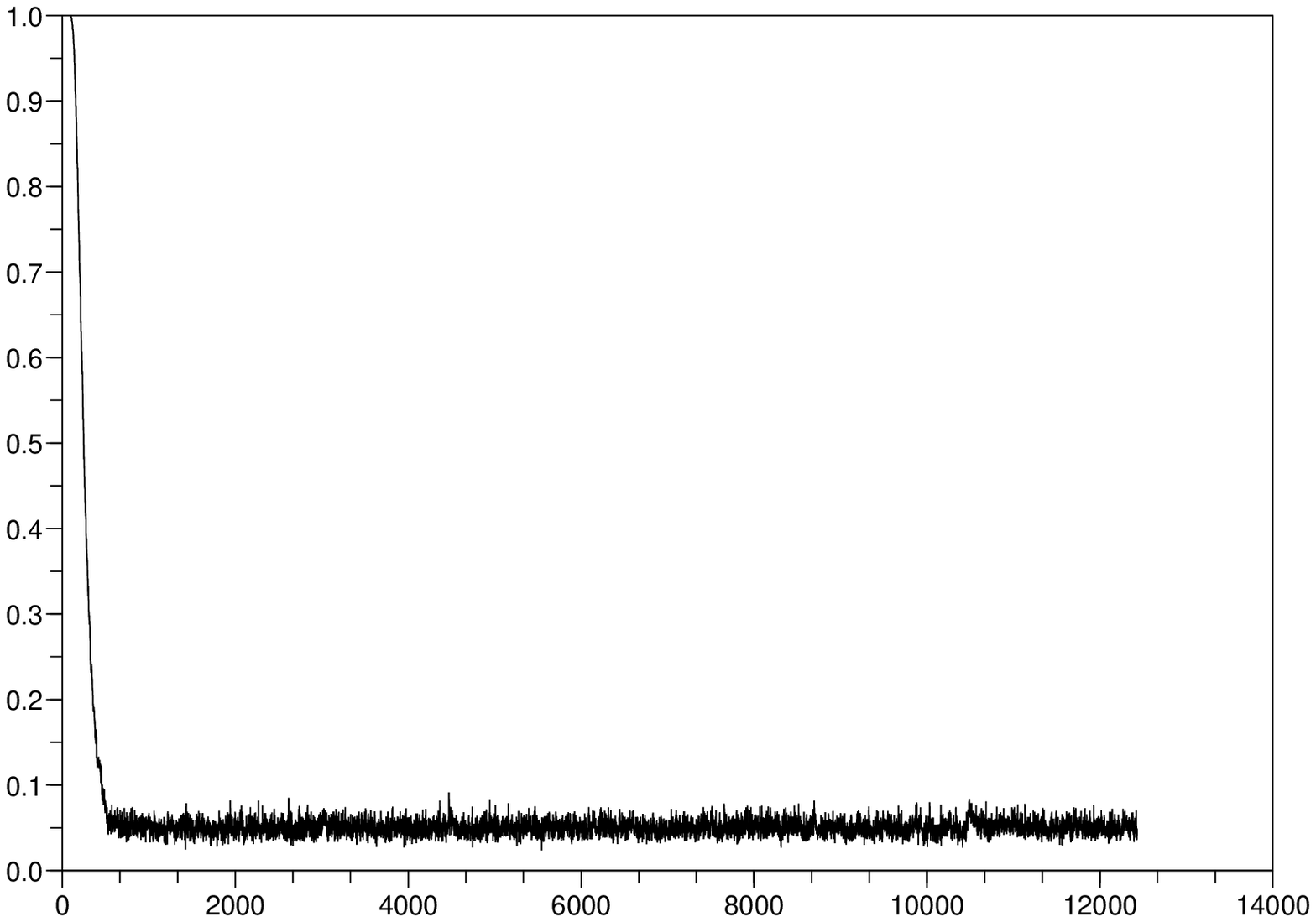}
  \caption{Probability distribution function of $\bar\phi(\nu^N)$
    and temporal evolution of $\bar d_0^N(t)$
    for case 2 and $\beta = 4$.}
  \label{fig4.2}
\end{figure}

\begin{figure}[H]
  \includegraphics[width=0.48\textwidth]{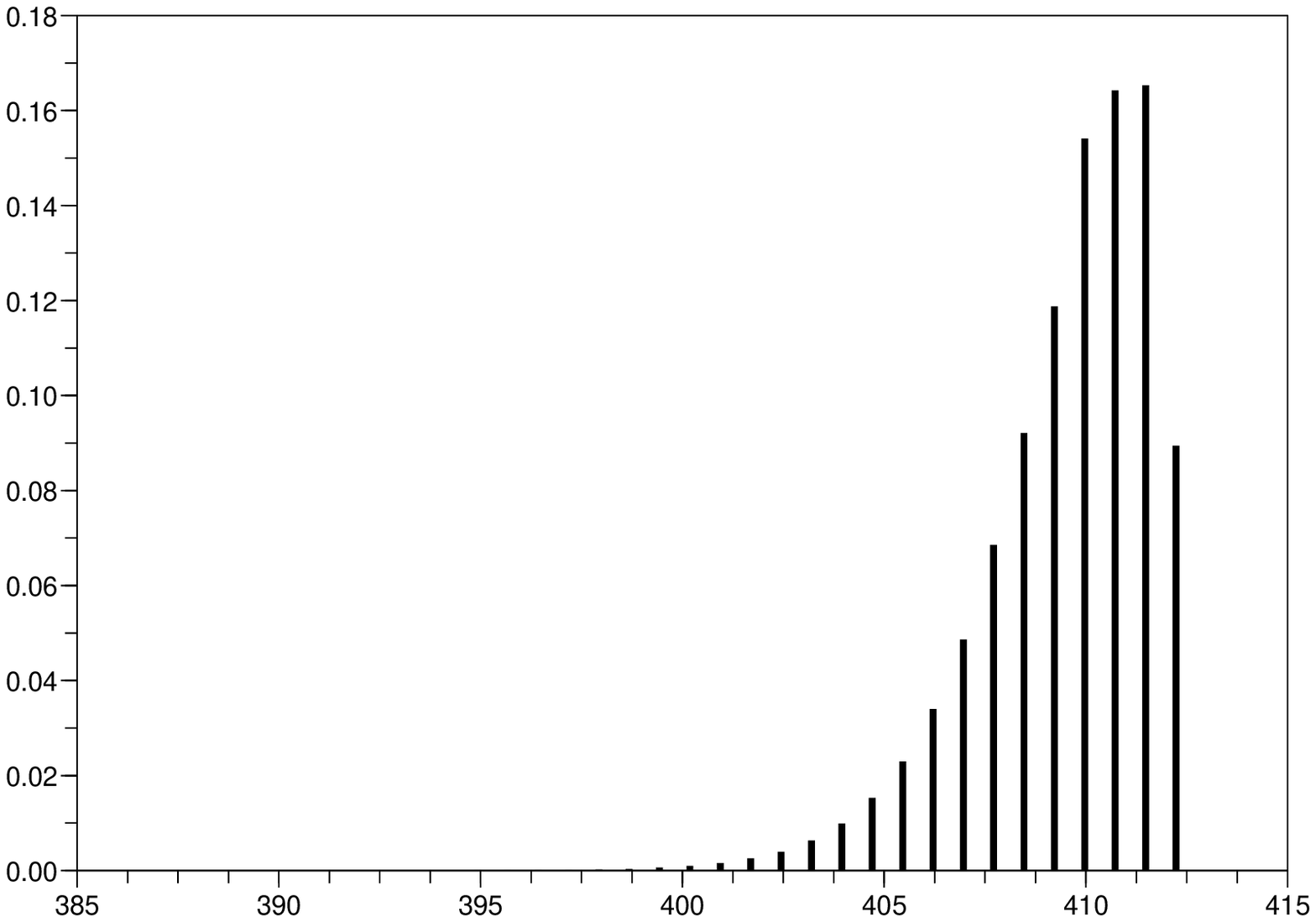}
  \includegraphics[width=0.48\textwidth]{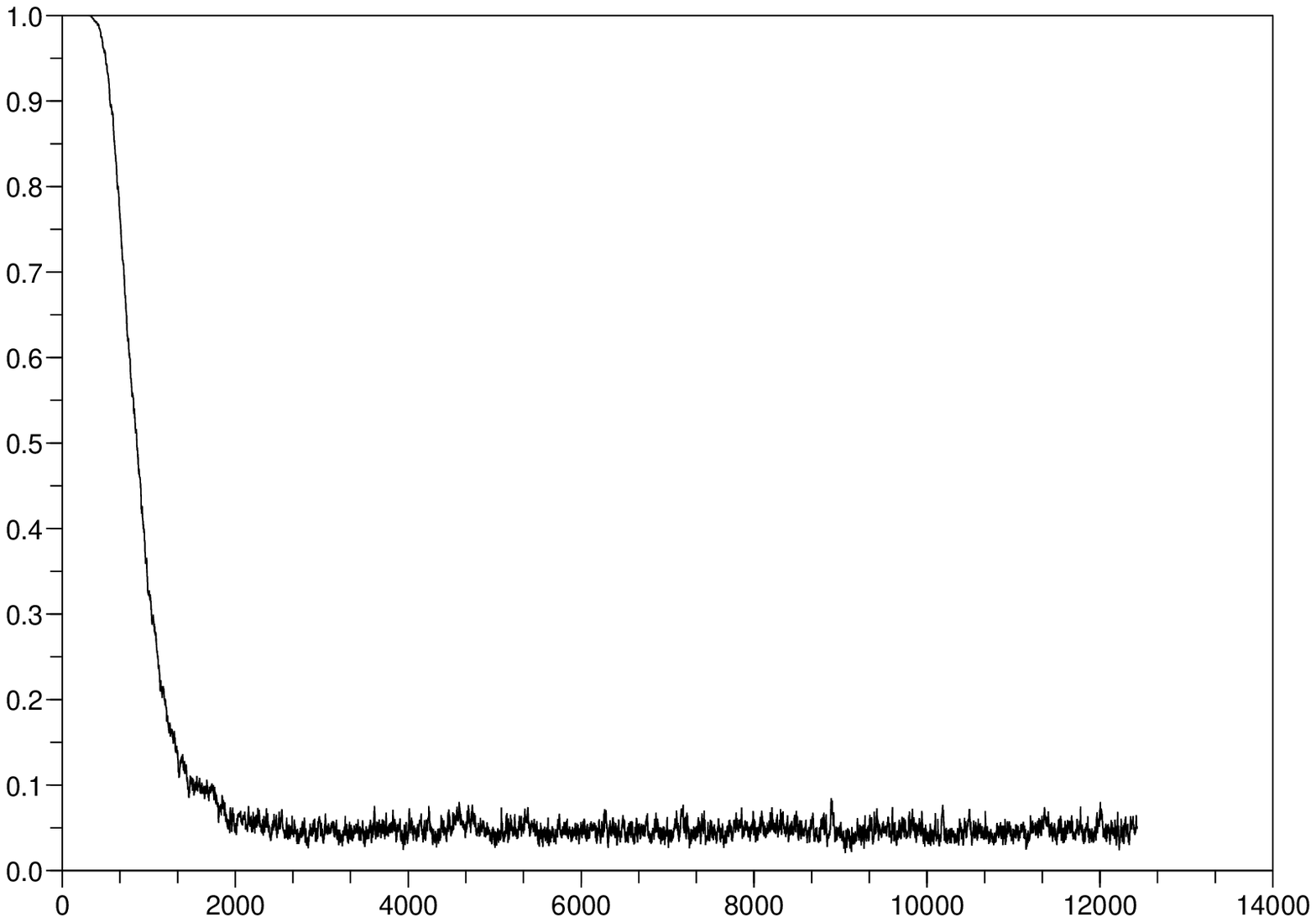}
  \caption{Probability distribution function of $\bar\phi(\nu^N)$
    and temporal evolution of $\bar d_0^N(t)$ 
    for case 2 and $\beta = 16$.}
  \label{fig4.3}
\end{figure}

\begin{figure}[H]
  \includegraphics[width=0.48\textwidth]{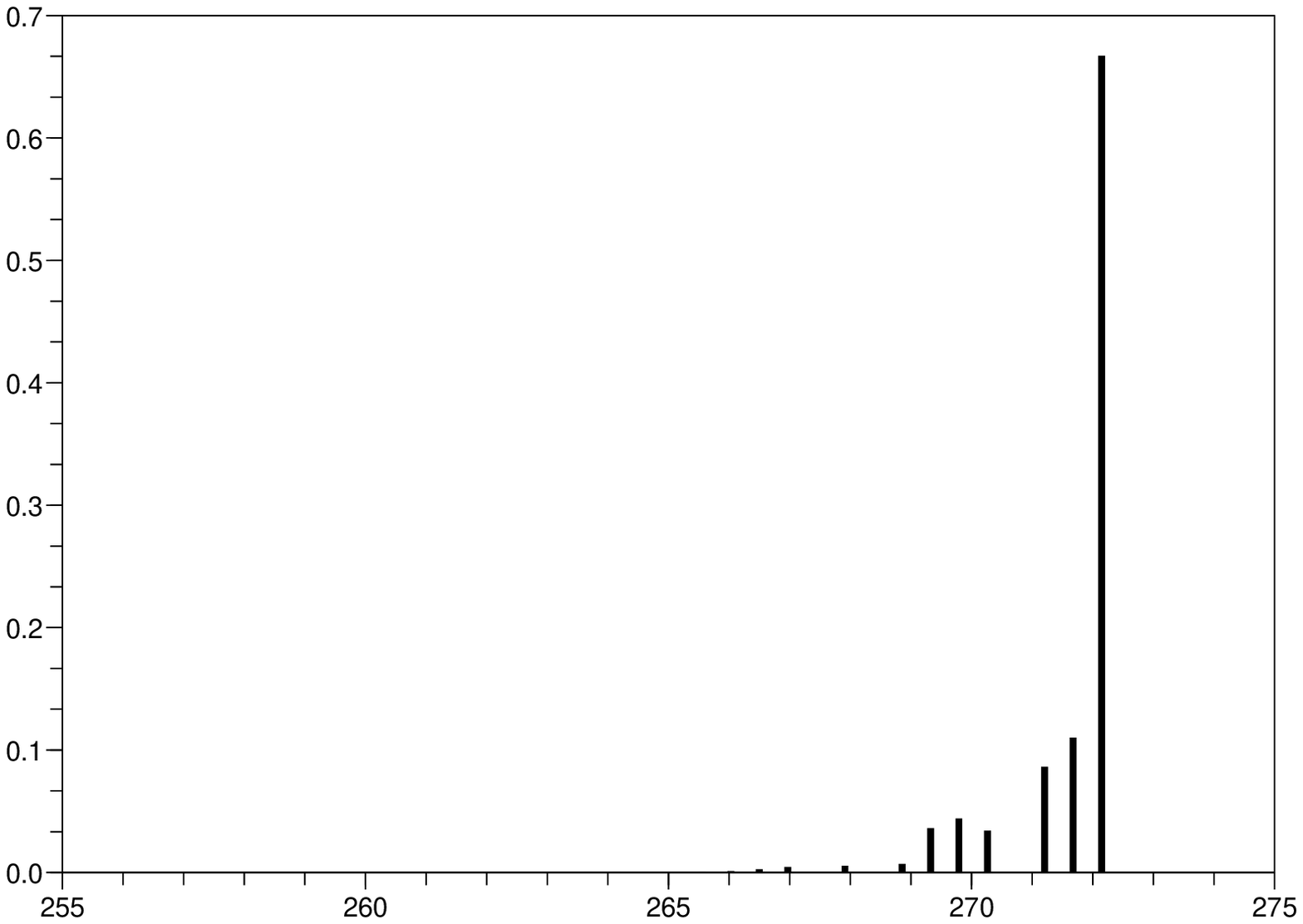}
  \includegraphics[width=0.48\textwidth]{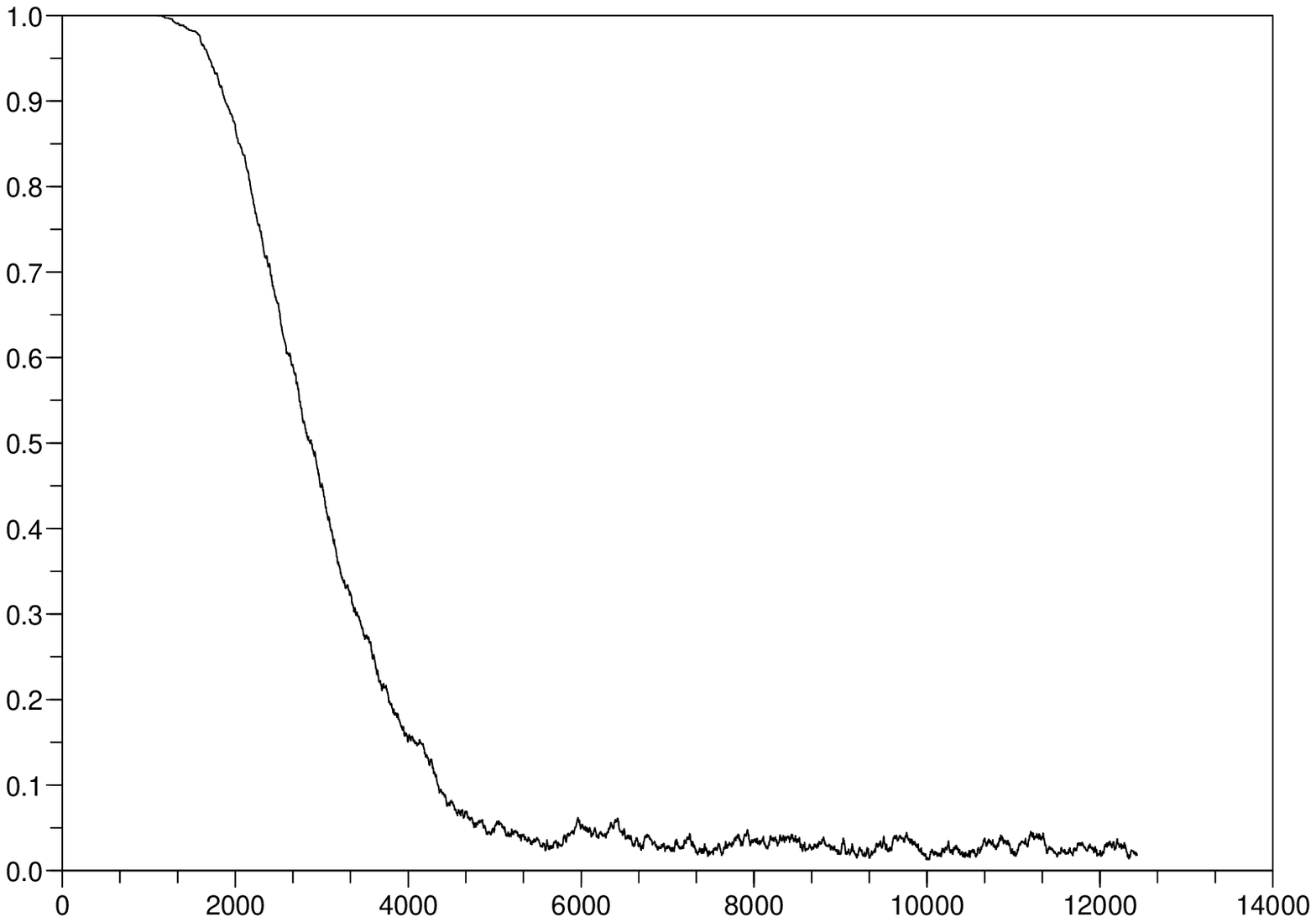}
  \caption{Probability distribution function of $\bar\phi(\nu^N)$
    and temporal evolution of $\bar d_0^N(t)$
    for case 2 and $\beta = 63$.}
  \label{fig4.4}
\end{figure}

\begin{figure}[H]
  \includegraphics[width=0.48\textwidth]{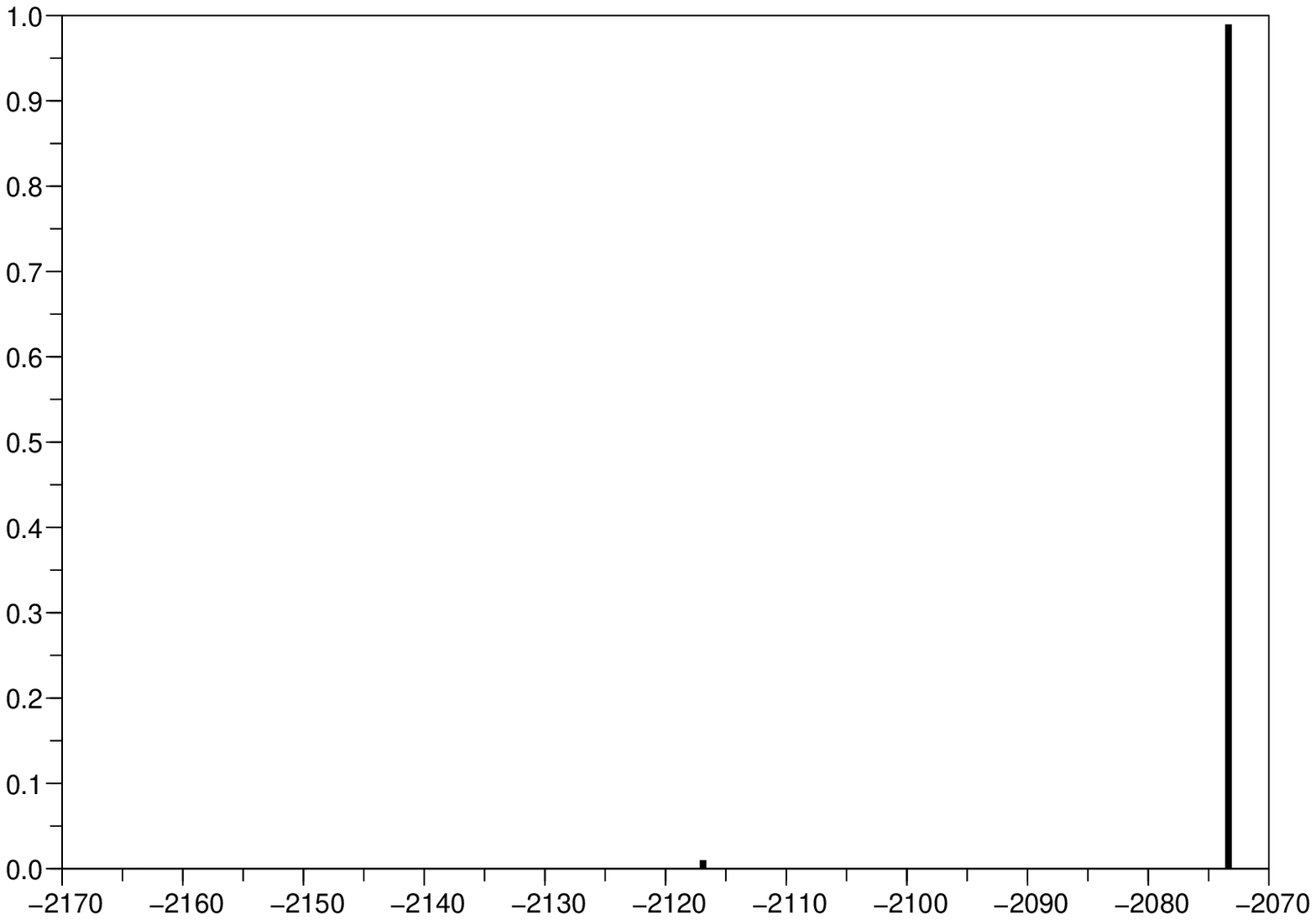}
  \includegraphics[width=0.48\textwidth]{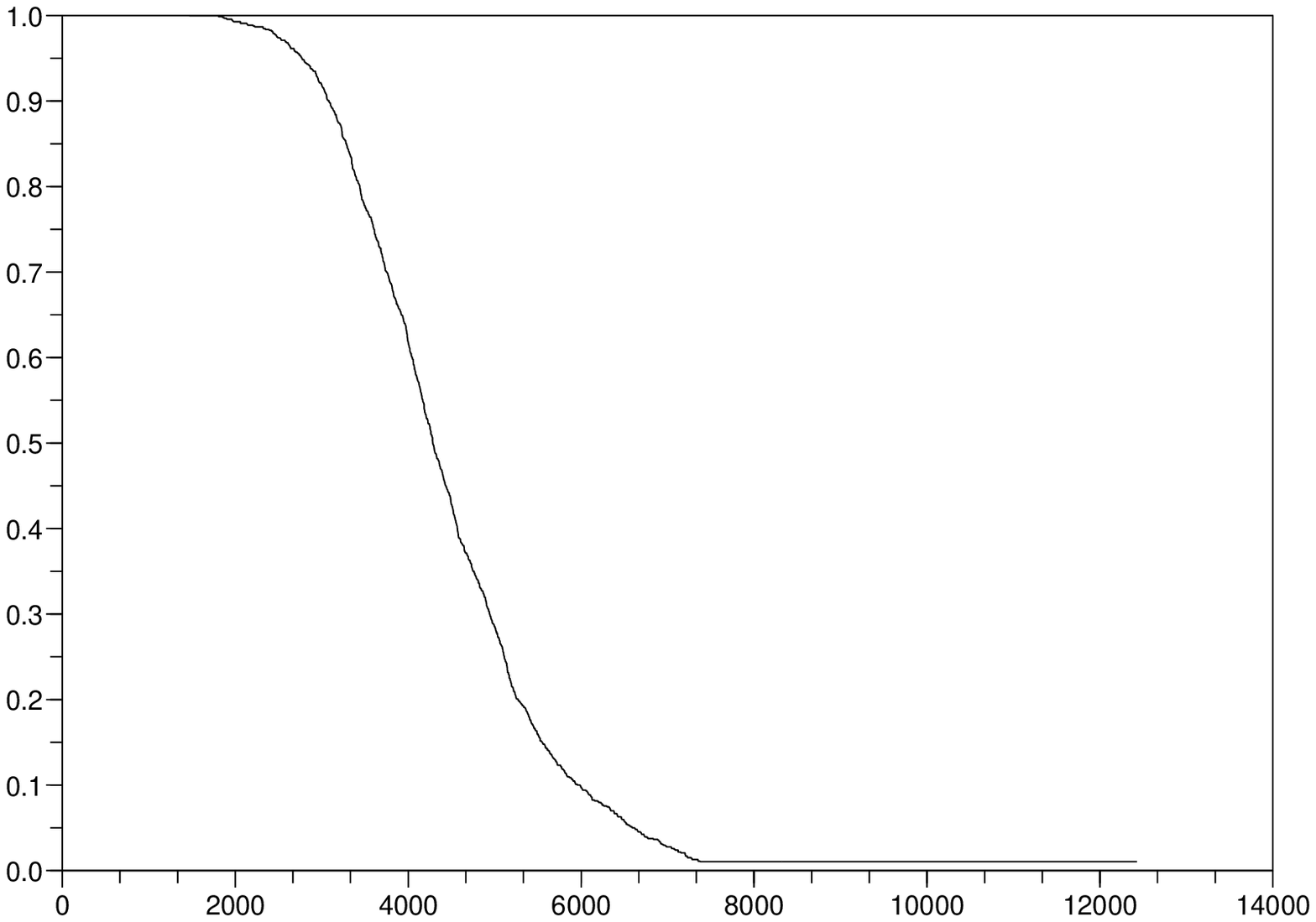}
  \caption{Probability distribution function of $\bar\phi(\nu^N)$
    and temporal evolution of $\bar d_0^N(t)$
    for case 2 and $\beta = 1000$.}
  \label{fig4.5}
\end{figure}

Summing up, we approximated $d(t)$ by
\begin{equation}
\bar d_0^N(t) 
= \max_{\eta_0\in{\cal X}_{k,m}^0} 
\|\bar\phi(\mu_t^N) - \bar\phi(\nu^N)\|_{TV}
\end{equation}
where ${\cal X}_{k,m}^0$ is the two configuration set
described above.
We chose $N = 1024$ and $M = \sqrt{N} = 32$
and we plotted for different temperatures
a graphical representation of the law of 
$\bar\phi(\nu^N)$ and the temporal evolution
of $\bar d_0^N(t)$ for $0\leq t \leq 2T$.
See Figures~\ref{fig3.1} to~\ref{fig3.5} for case 1
and Figures~\ref{fig4.1} to~\ref{fig4.5} for case 2.

\begin{figure}[htbp]
  \includegraphics[width=0.48\textwidth]{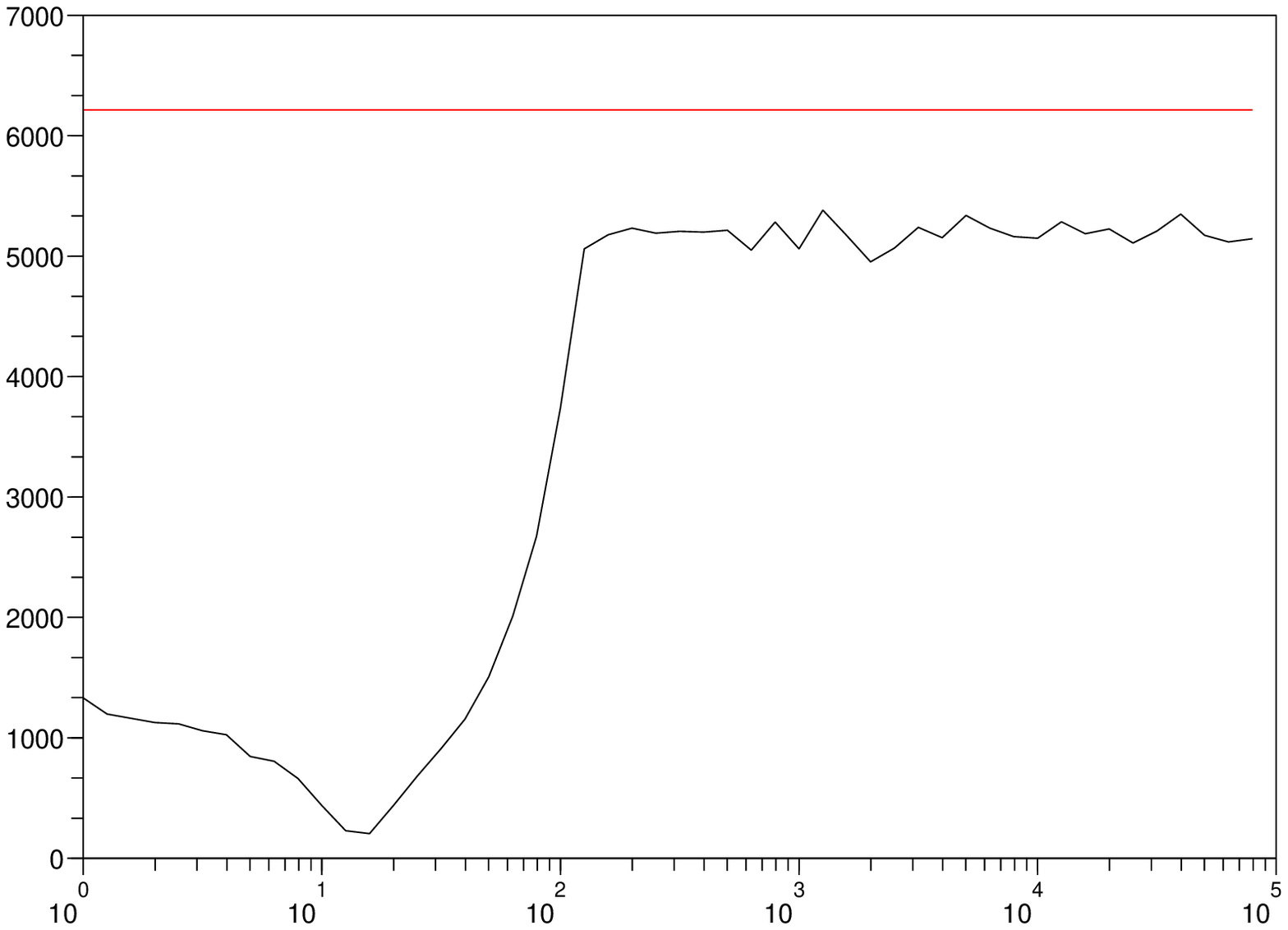}
  \includegraphics[width=0.48\textwidth]{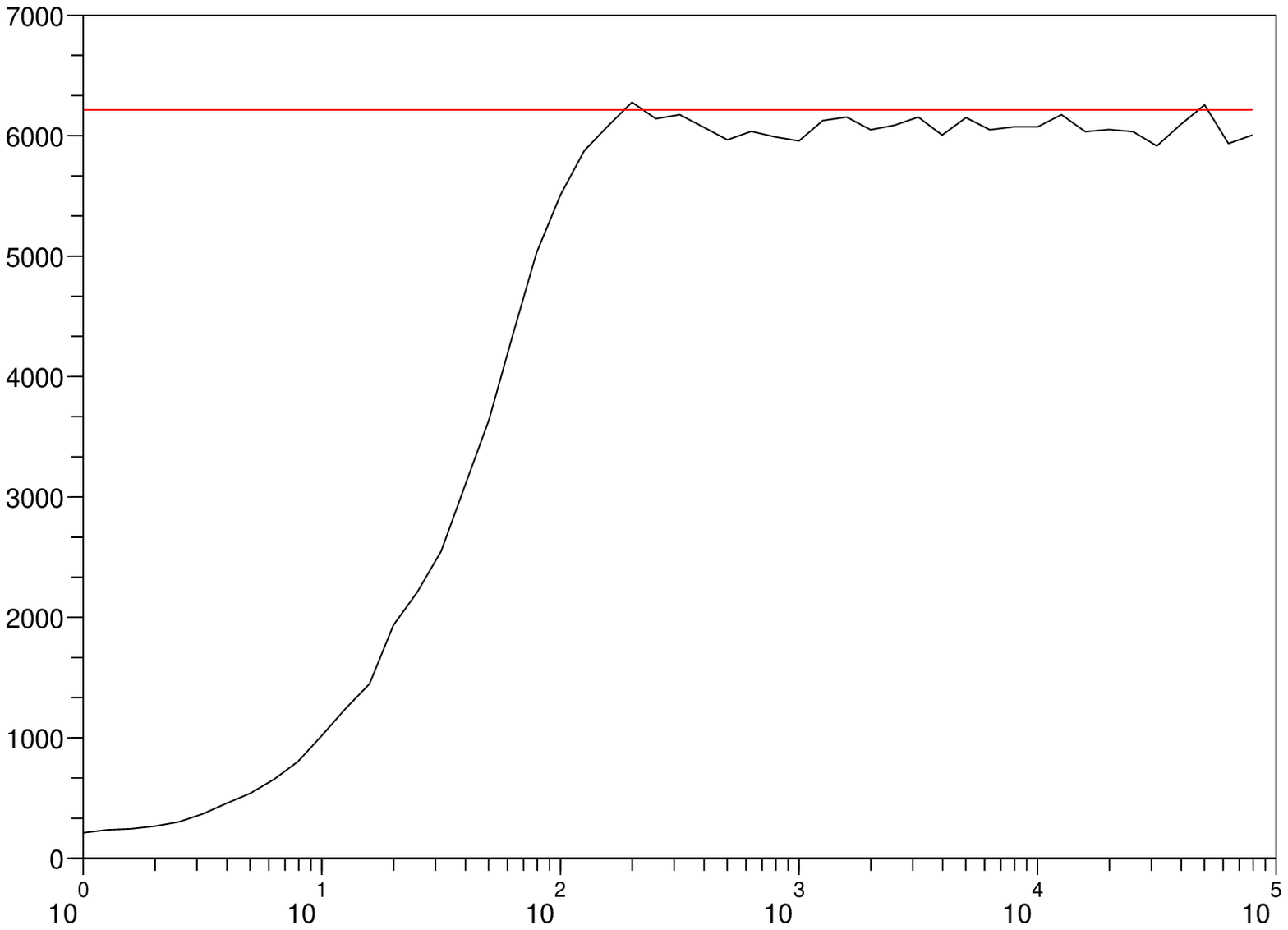}
  \caption{Dependency of $\hat t_\epsilon$ on $\beta$
    in case 1 and 2 respectively (semi-logarithmic representation).
    The horizontal line stands for our theoretical upper bound.} 
  \label{fig5}
\end{figure}

For each temperature we get in each case an estimation $\hat t_\epsilon$
of $t_\epsilon$. Estimating $t_\epsilon$ in this way for 
{51} temperatures ($\beta = 0$ included) we show 
in Figure~\ref{fig5} the graphical
representations of the simulated dependency of $\hat t_\epsilon$ on $\beta>0$
together with our upper bound from Theorem~\ref{thm1}.

We can now conclude with a few comments.
First, with such a procedure, $t_\epsilon$ 
is essentially underestimated. 
Indeed, except for the replacement
of $p^t(\eta_0,\cdot)$
and $\nu$ by $\mu^N_t$ and $\nu^N$
our successive approximations underestimated
$d(t)$ then $t_\epsilon$.
As a consequence we recover with case~2
that $km\ln(k/\epsilon)$
is the best upper bound on $t_\epsilon$
that is uniform over the disorder.

Second, we note that $\beta$ can be seen as
a tuning parameter for the disorder.
We recover that our theoretical estimate on $t_\epsilon$
can be improved by a factor of order $m$
in weak disorder situations.

Last, it is worth to note that case~1
suggests that $t_\epsilon$ is not
generally decreasing with the temperature.
In our preliminary simulations
we only considered $\eta_0 = (0, \dots, 0, 50)$
to estimate the mixing time.
Then we feared that our non monotonic
estimation could be an artefact of our simulation.
This led us to introduce the alternative saturated low level
starting configuration, confirming in that way
the non-monotonicity.
\end{appendix}

\clearpage
\section*{Acknowledgment}
J.R. is very grateful to everyone who made his stay in Rome possible and pleasant.

\end{document}